\documentclass[10pt,twocolumn,twoside]{IEEEtran}
\usepackage{textcomp}
\usepackage{amsmath,amssymb,pstricks,color,dsfont}
\usepackage[sort,numbers]{natbib}
\usepackage{psfrag,graphicx,epsfig}
\usepackage{algorithm,algorithmic,xspace}
\usepackage{subfigure}
\usepackage{float}
\usepackage{epstopdf}
\usepackage{savesym}
\usepackage{gensymb}

\usepackage[algo2e,vlined,boxruled,linesnumbered]{algorithm2e}

\allowdisplaybreaks

\newcommand{\bremark}{\begin{remark} \begin{rm} }
\newcommand{\eremark}{ \end{rm} \hfill\rule{2mm}{2mm}
\end{remark} }
\newcommand{\btheorem}{\begin{theorem} \begin{rm} }
\newcommand{\etheorem}{ \end{rm} \hfill\rule{2mm}{2mm}
\end{theorem} }
\newcommand{\blemma}{\begin{lemma} \begin{rm} }
\newcommand{\elemma}{ \end{rm} \hfill\rule{2mm}{2mm}
\end{lemma} }
\newcommand{\bcorollary}{\begin{corollary} \begin{rm} }
\newcommand{\ecorollary}{ \end{rm} \hfill\rule{2mm}{2mm}
\end{corollary} }
\newcommand{\bdefinition}{\begin{definition}\begin{rm} }
\newcommand{\edefinition}{ \end{rm} \hfill\rule{2mm}{2mm}
\end{definition} }
\newcommand{\bproposition}{\begin{proposition} \begin{rm} }
\newcommand{\eproposition}{ \end{rm} \hfill\rule{2mm}{2mm}
\end{proposition} }
\newcommand{\bexample}{\begin{example} \begin{rm} }
\newcommand{\eexample}{ \end{rm} \hfill\rule{2mm}{2mm}
\end{example} }
\newcommand{\basm}{\begin{assumption} \begin{rm}}
\newcommand{\easm}{\end{rm} \hfill\rule{2mm}{2mm}
\end{assumption}}
\newcommand{\bclaim}{\begin{claim} \begin{rm} }
\newcommand{\eclaim}{ \end{rm} \hfill\rule{2mm}{2mm}
\end{claim} }

\newtheorem{claim}{\bf Claim}[section]
\newtheorem{theorem}{\bf Theorem}[section]
\newtheorem{lemma}{\bf Lemma}[section]
\newtheorem{definition}{\bf Definition}[section]
\newtheorem{remark}{\bf Remark}[section]
\newtheorem{corollary}{\bf Corollary}[section]
\newtheorem{proposition}{\bf Proposition}[section]
\newtheorem{example}{\bf Example}[section]
\newtheorem{assumption}{\bf Assumption}[section]
\newcommand\oprocendsymbol{\hbox{$\bullet$}}
\newcommand\oprocend{\relax\ifmmode\else\unskip\hfill\fi\oprocendsymbol}

\pdfoutput=1

\begin{document}

\title{Distributed economic control of dynamically coupled networks}

\author{Yang~Lu~and~Minghui~Zhu
\thanks{Y. Lu and M. Zhu are with the School of Electrical Engineering and Computer Science, Pennsylvania State University, 201 Old Main, University Park, PA, 16802, {\tt\small yml5046@psu.edu, muz16@psu.edu}. This work was partially supported by NSF grant ECCS-1710859.}}

\maketitle

\begin{abstract}
This paper investigates the synthesis of distributed economic control algorithms under which dynamically coupled physical systems are regulated to a variational equilibrium of a constrained convex game. We study two complementary cases: (i) each subsystem is linear and controllable; and (ii) each subsystem is nonlinear and in the strict-feedback form. The convergence of the proposed algorithms is guaranteed using Lyapunov analysis. Their performance is verified by two case studies on a multi-zone building temperature regulation problem and an optimal power flow problem, respectively.
\end{abstract}


\section{Introduction}

Control systems have been experiencing dramatic changes where advanced information and communication technologies are integrated to improve system performance and produce new functions. This class of new-generation control systems is referred to as cyber-physical systems (CPSs), e.g., mobile robotic networks, the smart grid, smart buildings and intelligent transportation systems. The entities of CPSs are spatially distributed and thus distributed coordination becomes necessary to achieve network-wide goals \cite{fb-jc-sm:09,MZ-SM-distributed-book:2015,MM-ME:10}. 

Currently, many large-scale CPSs, e.g., power systems, adopt top-down hierarchical control. The decision-making layer first determines an optimal operating point via solving, e.g., an optimization problem. The physical layer waits for the convergence of the decision-making layer and then regulates physical dynamics to the desired operating point. For example, in optimal power flow (OPF) problems \cite{MH-DX-LW:2016,JZ:2015} 
the decision-making layer solves an economic dispatch problem to schedule the nominal operating levels of active generators and inter-area power flows. Once the economic dispatch problem is solved, the physical layer performs frequency control to adjust the set-points of the generators to drive their frequencies back to their nominal values and drive the inter-area power flows to the scheduled values.

One potential issue of the hierarchical approach is that, in an iterative seeking of an optimal operating point, the decision-making layer needs to evaluate quality of intermediate estimates so as to determine whether the seeking process terminates or not.
The quality evaluation could be challenging, especially in a distributed setting, in which individual agents may lack of information to verify solution quality for a network-wide problem.

As pointed out in the literature of economic model predictive control (MPC) \cite{ME-JL-PDC:2017,JBR-DA-CNB:2012}, another potential issue of hierarchical control is efficiency loss. A general reason for the issue is that optimal operating points change more frequent and the time scales on which physical dynamics react to the changes shrink. The issue motivates recent studies of economic MPC as a complement of hierarchical control. In economic MPC, the controller directly optimizes in real time the economic performance of physical dynamics. However, the design of economic MPC still requires to decide the optimal steady state of physical dynamics in advance.

Alternatively, there have been recent efforts on the approach where the physical layer keeps track of intermediate results of the decision-making layer all the time. In this paper, the approach is referred to as economic control. Since the decision-making layer and the physical layer operate in a single time scale, solution evaluation is resolved and physical dynamics are more responsive to dynamic changes of optimal operating points.
In the smart grid, \cite{EM-CZ-SL:2017} proposed a
strategy for linear power systems which, by using virtual flows to constrain real flows without altering the primal-dual interpretation of the network dynamics, can simultaneously address frequency control and economic dispatch.
In multi-zone building systems, \cite{TH-XZ-WS-MZ-NL:2017CDC,TH-XZ-WS-MZ-NL:2017CCTA} proposed an integrated design such that the optimization of the decision-making layer is formulated as a dynamical system and the physical layer keeps track of the dynamic feedback from the decision-making layer.
The papers aforementioned study specific applications. On the other hand, recent papers \cite{Stankovic.Johansson.Stipanovic:11,PF-MK-TB:11} provide general theories to integrate physical dynamics into game-theoretic learning for dynamically decoupled networks, e.g., mobile robotic networks. In \cite{Stankovic.Johansson.Stipanovic:11}, the agents' actions are filtered through decoupled stable linear time-invariant filters
before affecting objective functions.
In \cite{PF-MK-TB:11}, the agents' actions act as inputs to a stable nonlinear dynamic system whose outputs are the agents' costs.
In contrast, this paper studies dynamically coupled networked systems which are unstable without controls and need to be stabilized in a distributed manner. Moreover, in \cite{Stankovic.Johansson.Stipanovic:11,PF-MK-TB:11}, no constraint was included in the games of interest.
This paper investigates a class of convex games, which were first introduced in \cite{Arrow.Debreu:54} and have received substantial attentions \cite{Facchinei.Kanzow:07}.
A number of algorithms have been proposed to compute Nash equilibria of convex games, including, to name a few, ODE-based methods \cite{Rosen:65}, nonlinear Gauss-Seidel-type approaches \cite{Pang.Scutari.Facchinei.Wang:08}, iterative primal-dual Tikhonov schemes \cite{Yin.Shanbhag.Mehta:11}, consensus-based methods \cite{JK-AN-UVS:2016}, best-response dynamics \cite{PE:10} and finite difference approximation \cite{MZ-EF:Auto14}. 
The learning algorithms of game theory have been applied to, e.g., networking \cite{Altman.Basar:98}, transportation \cite{LB-MC-SG:2016}, multi-robot systems \cite{Zhu.Martinez:12}, power systems \cite{WS-ZH-HVP-TB:2012} and multi-zone building systems \cite{SAA-WHS-GZ:2018}. In all these papers, the game-theoretic learning algorithms ignore real-time physical dynamics which feature CPSs.


{\color{blue}
}

\emph{Contribution.} In this paper, we consider the synthesis of distributed economic control algorithms under which dynamically coupled physical systems are regulated to a variational equilibrium (VE), a special type of Nash equilibrium (NE), of a constrained convex game. The game includes distributed separable optimization as a special case. We investigate two complementary cases where the physical dynamics are linear and nonlinear, respectively.

For linear dynamics with controllable subsystems, we propose an indirect approach for the distributed economic control algorithm design. We synthesize distributed update rules for auxiliary variables such that they asymptotically converge to a VE of the game of interest, but their evolution is not restricted by the inherent physical dynamics. The distributed decision-making update rule extends the projected discontinuous dynamics in \cite{AC-EM-JC:2016} from convex optimization to convex games. We then design a novel distributed control law such that the states and inputs of the physical dynamics keep track of the auxiliary variables and thus asymptotically converge to the same VE of the game. The designed algorithm is verified by a multi-zone building temperature regulation problem. The simulation demonstrates that, compared with the hierarchical control scheme, the proposed economic control scheme can significantly reduce efficiency loss.

For nonlinear dynamics in the strict-feedback form,
we develop a direct approach without introducing auxiliary variables. In particular, we perform a sequence of coordinate transformations such that each subsystem is converted into the one where its first state is driven by a primal-dual dynamics and it is followed by a chain of integrators. We then design a distributed control law such that the integrators are stabilized and the primal-dual dynamics asymptotically converges to a VE of the game of interest. The algorithm and its analysis developed for this case are novel. The proposed algorithm is verified by an optimal power flow problem.

Preliminary results of this paper were published in \cite{YL-MZ:NECSYS2015Privacy}, where the case of nonlinear dynamics was not investigated, no case study was included, and the proofs were omitted.

\emph{Notations and notions.}
Let $(x_i)_{i=1}^N=[x_1^T,\cdots,x_N^T]^T$ and $x_{-i}=[x_1^T,\cdots,x_{i-1}^T,x_{i+1}^T,\cdots,x_N^T]^T$. Let $\|  \cdot  \|_1$ and $\|  \cdot  \|$ denote the 1-norm and 2-norm of a vector or a matrix, respectively. Denote the zero column vector of size $n$ by $0_n$ and the identity matrix of size $n\times n$ by $I_n$. Let $\mathbb{R}_+^n$ denote the set of nonnegative real column vectors of size $n$.
Denote by $\mathrm{diag}(A_i)_{i=1}^M$ the block-diagonal matrix composed by sub-matrices $A_1,\cdots,A_M$, such that the $j$-th diagonal block is $A_j$ and all the off-diagonal blocks are zero matrices. Given a function $f:\mathbb{R}^n\to\mathbb{R}$ of $x=(x_i)_{i=1}^n\in\mathbb{R}^n$, for a nonnegative integer $\ell$, we say that $f$ is of class $C^\ell$ if for all possible combinations of $\ell_1,\ell_2,\cdots,\ell_n$, where each of $\ell_1,\ell_2,\cdots,\ell_n$ is an integer between $0$ and $\ell$ such that $\ell_1+\ell_2+\cdots+\ell_n=\ell$, the partial derivative $\frac{\partial^\ell f}{\partial x_{1}^{\ell_1}\partial x_{2}^{\ell_2}\cdots\partial x_{n}^{\ell_n}}$ exists and is continuous. Given a function $V:\mathbb{R}^n\to\mathbb{R}$ and $\delta>0$, denote the sublevel set of $V$ by $V^{-1}(\leq\delta)=\{x\in\mathbb{R}^n:V(x)\leq\delta\}$.

Let $K\subseteq\mathbb{R}^n$ be a closed and convex set. Denote the boundary, interior and closure of $K$ by ${\rm bd}(K)$, ${\rm int}(K)$ and ${\rm cl}(K)$, respectively. Given $x\in{\rm bd}(K)$, define the normal cone of $K$ at $x$
by $N_K(x)=\{\gamma:\|\gamma\|=1,\;{\rm and}\;\gamma^T(x-y)\leq0,\,\forall y\in K\}$. Given $x\in\mathbb{R}^n$, the point projection of $x$ onto $K$ is ${\rm Proj}_K(x)={\rm argmin}_{z\in K}\|z-x\|$. Given $x\in K$ and $v\in\mathbb{R}^n$, the vector projection of $v$ at $x$ with respect to $K$ is $\Pi_K(x,v)=\lim_{\delta\to0^+}\frac{{\rm Proj}_K(x+\delta v)-x}{\delta}$.
The following lemma is adopted from Lemma 2.1 of \cite{AN-DZ:1996}.

\blemma
\label{lemma vector projection}
If $x\in {\rm int}(K)$, then $\Pi_K(x,v)=v$; if $x\in{\rm bd}(K)$, then $\Pi_K(x,v)=v+\beta(x)\gamma^*(x)$, where $\gamma^*(x)={\rm argmin}_{\gamma\in N_K(x)}v^T\gamma$ and $\beta(x)=\max\{0,-v^T\gamma^*(x)\}$; and $\|\Pi_K(x,v)\|\leq\|v\|$ for any $x\in K$ and any $v\in\mathbb{R}^n$.
\elemma

We next review some necessary concepts on projected discontinuous dynamical systems. Given a set $K\subseteq\mathbb{R}^n$ and a map $F:K\to\mathbb{R}^n$, consider the projected differential equation
\begin{align}
\label{differential inclusion}
\dot x=\Pi_K(x,F(x)).
\end{align}

Given $x^0\in K$, a solution of \eqref{differential inclusion} on the interval $[0,T)\subseteq\mathbb{R}$ is a map $\Gamma:[0,T)\to K$ with $\Gamma(0)=x^0$ that is absolutely continuous on $[0,T)$ and satisfies $\dot \Gamma(t)=\Pi_K(\Gamma(t),F(\Gamma(t)))$ almost everywhere on $[0,T)$.
%
%
A set $S\subseteq\mathbb{R}^n$ is invariant under \eqref{differential inclusion} if every solution starting from any point in $S$ remains in $S$. For a solution $\Gamma$ of \eqref{differential inclusion} defined on $[0,\infty)$, the omega-limit set $\Omega(\Gamma)$ is defined by $\Omega(\Gamma)=\{y\in\mathbb{R}^n:\exists\{t_k\}_{k=1}^\infty\subseteq[0,\infty)\;{\rm with}\;\lim_{k\to\infty}t_k=\infty\;{\rm and}\;\lim_{k\to\infty}\Gamma(t_k)=y\}$. Given a continuously differentiable function $V:K\to\mathbb{R}$, the Lie derivative of $V$ along \eqref{differential inclusion} at $x\in K$ is $\mathcal{L}_{K,F}V(x)=\nabla V(x)^T\Pi_K(x,F(x))$. The next result adopted from \cite{AB-FC:2006} provides an invariance principle for \eqref{differential inclusion}.

\blemma
\label{lemma invariance}
Let $S\subseteq\mathbb{R}^n$ be compact and invariant under \eqref{differential inclusion}. Assume that for any $x^0\in S$ there exists a unique solution of \eqref{differential inclusion} starting at $x^0$ and that its omega-limit set is invariant under \eqref{differential inclusion}. Suppose that there exists a continuously differentiable function $V:K\to\mathbb{R}$ such that $V$ is positive definite and $\mathcal{L}_{K,F}V(x)\leq0$ for all $x\in S$. Then any solution of \eqref{differential inclusion} starting from $S$ converges to the largest invariant set contained in ${\rm cl}(\{x\in S:\mathcal{L}_{K,F}V(x)=0\})$.
\elemma

We adopt the following notions of strictly and strongly monotone maps from \cite{FF-JSP:03}.

\bdefinition
A map $\Phi:\mathbb{R}^r\to\mathbb{R}^r$ is strictly monotone on $W\subseteq\mathbb{R}^r$ if $(w-w')^T(\Phi(w)-\Phi(w'))>0$ for all $w,w'\in W$ with $w\neq w'$. The map $\Phi$ is strongly monotone on $W$ if there exists $M>0$ such that $(w-w')^T(\Phi(w)-\Phi(w'))\geq M\|w-w'\|^2$ for all $w, w'\in W$.
\label{def: strict monotone}
\edefinition

\section{Motivating examples\label{section motivation}}

In this section, we use a multi-zone building temperature regulation problem and an optimal power flow (OPF) problem to motivate our problem formulation.

\subsection{Multi-zone building temperature regulation\label{section MZB}}

Consider a set of zones $\mathcal{V}=\{1,\cdots,N\}$ and an air handling unit (AHU). The meanings of the parameters and variables are listed in Table \ref{thermal parameters}. For each $i\in\mathcal{V}$, let $\mathcal{N}_i^{PH}\subseteq\mathcal{V}$ be the set of neighboring zones of zone $i$. For each zone $i$, the following dynamic model is adopted from \cite{YM-GA-FB:2011}:
\begin{align}
\label{thermal dynamics}
&C_{i1}\dot T_{i1}(t)=m_i^s(t)c_p(T_i^s(t)-T_{i1}(t))+\frac{T_{i2}(t)-T_{i1}(t)}{R_i}\nonumber\\
&+\sum\nolimits_{j\in\mathcal{N}_i^{PH}}\frac{T_{j1}(t)-T_{i1}(t)}{R_{ji}}+\frac{T^{oa}(t)-T_{i1}(t)}{R_{i}^{oa}}+P_i^d(t),\nonumber\\
&C_{i2}\dot T_{i2}(t)=\frac{T_{i1}(t)-T_{i2}(t)}{R_{i}},\nonumber\\
&T_i^s(t)=\delta(t)\frac{\sum_{j\in\mathcal{V}}m_j^s(t)T_{j1}(t)}{\sum_{j\in\mathcal{V}}m_j^s(t)}+(1-\delta(t))T^{oa}(t)\nonumber\\
&-\Delta T^c(t)+\Delta T_i^h(t).
\end{align}

\begin{table}[H]
\renewcommand{\arraystretch}{1.1}
\caption{Parameters/variables of the multi-zone building problem}
\label{thermal parameters}
\centering
\begin{tabular}{|l|l|}
\hline
$C_{i1}/C_{i2}$ & thermal capacitance of fast/slow-dynamic mass at zone $i$\\
\hline
$T_{i1}/T_{i2}$ & temperature of fast/slow dynamic mass at zone $i$\\
\hline
$T_i^r$ & reference temperature of zone $i$\\
\hline
$T_i^s$ & temperature of air supplied to zone $i$\\
\hline
$m_i^s$ & mass flow rate of air supplied to zone $i$\\
\hline
$T^{oa}$ & outside air temperature\\
\hline
$\Delta T^c$ & temperature difference across AHU cooling coil\\
\hline
$\Delta T_i^h$ & temperature difference across reheat coil at zone $i$\\
\hline
$\delta$ & AHU recirculation damper position\\
\hline
$P_i^d$ & external load\\
\hline
$R_i$ & heat resistance between $C_{i1}$ and $C_{i2}$\\
\hline
$R_{ji}$ & thermal resistance between zone $i$ and zone $j$\\
\hline
$R_i^{oa}$ & thermal resistance between zone $i$ and outside air\\
\hline
$c_p$ & specific heat capacity of zone air\\
\hline
$\underline T_i$/$\overline T_i$ & lower/upper limit of the temperature of zone $i$\\
\hline
$\underline{\Delta T}_i^h$/$\overline{\Delta T}_i^h$ & lower/upper limit of $\Delta T_i^h$\\
\hline
\end{tabular}
\end{table}
\normalsize

The agents aim to minimize the deviations between the zone temperatures and their references with minimum control efforts. Meanwhile, the zone temperatures and controls are subject to some hard limits and the coupled constraints on thermal balance. All these objectives are encapsulated in the following optimization problem which is adopted from \cite{TH-XZ-WS-MZ-NL:2017CDC}:
\begin{align}
\label{thermal OP}
&\!\!\!\!\min\limits_{T_{1},T_{2},\Delta T^h}\sum\limits_{i\in\mathcal{V}}(c_{i1}^x(T_{i1}-T_i^{r})^2+c_{i2}^x(T_{i2}-T_i^{r})^2+c_{i}^u(\Delta T_i^h)^2)\nonumber\\
&\!{\rm s.t.}\;\; m_i^sc_p(T_i^s-T_{i1})+\frac{T_{i2}-T_{i1}}{R_i}+\sum\nolimits_{j\in\mathcal{N}_i^{PH}}\frac{T_{j1}-T_{i1}}{R_{ji}}\nonumber\\
&\quad\;\;+\frac{T^{oa}-T_{i1}}{R_{i}^{oa}}+P_i^d=0,\;\forall i\in\mathcal{V}\nonumber\\
&\quad\;\; T_{i1}-T_{i2}=0,\;\forall i\in\mathcal{V},\nonumber\\
&\quad\;\;\underline{\Delta T}_i^h\leq \Delta T_i^h\leq \overline{\Delta T}_i^h,\;\underline T_i\leq T_{i1},T_{i2}\leq \overline T_i,\;\forall i\in\mathcal{V}
\end{align}
where $c_{i1}^x$, $c_{i2}^x$ and $c_{i}^u$ are constant positive weights associated with $T_{i1}$, $T_{i2}$ and $\Delta T_i^h$, respectively.

In traditional hierarchical control, the decision-making layer first solves problem \eqref{thermal OP} to obtain an optimal operating point $(T_{1}^*,T_{2}^*,\Delta T^{h*})$. After the decision-making layer converges, the physical layer performs temperature regulation to steer the physical dynamics \eqref{thermal dynamics} to the desired operating point \cite{FAQ-CNJ:2018}. To address the issues of hierarchical control, we aim to design an economic control algorithm such that the temperature regulation in the physical dynamics \eqref{thermal dynamics} could keep track of the intermediate solutions of problem \eqref{thermal OP}.


\subsection{Optimal power flow\label{section OPF}}

Consider a power network comprising a set of power generators $\mathcal{V}=\{1,\cdots,N\}$. The meanings of the parameters and variables are listed in Table \ref{OPF parameters}. For each $i\in\mathcal{V}$, let $\mathcal{N}_i^{PH}\subseteq\mathcal{V}$ identify the set of neighboring generators of generator $i$. For convenience of notation, for any pair $(i,j)$, let $\theta_{ij}(t)=\theta_i(t)-\theta_j(t)$. The following inherent physical dynamics of generator $i$ is adopted from \cite{YW-DJH-GG:1998}:
\begin{align}
\label{motivating example dynamics}
&\dot\theta_i(t)=\omega_i(t)-\tilde\omega,\nonumber\\
&\dot\omega_i(t)=-\frac{D_i}{2T_{H_i}}\omega_i(t)+\frac{\omega_0}{2T_{H_i}}(P_{m_i}(t)-P_{m_i0})\nonumber\\
&\qquad\quad\;-\sum\nolimits_{j\in\mathcal{N}_i^{PH}}\frac{\omega_0t_{ij}}{2T_{H_i}}(\sin\theta_{ij}(t)-\sin\theta_{ij}(0)),\nonumber\\
&\dot P_{m_i}(t)=-\frac{1}{T_{M_i}}(P_{m_i}(t)-P_{m_i0})+\frac{K_{M_i}}{T_{M_i}}(P_{E_i}(t)-P_{E_i0}),\nonumber\\
&\dot P_{E_i}(t)=-\frac{K_{E_i}\omega_i(t)}{T_{E_i}R_i\omega_0}-\frac{1}{T_{E_i}}(P_{E_i}(t)-P_{E_i0}-u_i(t))
\end{align}
where $\tilde\omega=50Hz$ is the desirable relative speed.

\begin{table}[!t]
\renewcommand{\arraystretch}{1.1}
\caption{Parameters/variables of the OPF problem}
\label{OPF parameters}
\centering
\begin{tabular}{|l|l||l|l|}
\hline
$a,b,c$ & cost parameters & $P_C$ & power control input\\
\hline
$T_{H}$ & inertia constant & $u_i$ & $P_{C_i}-P_{m_i0}$\\
\hline
$D$ & damping coefficient & $\omega$ & relative speed\\
\hline
$T_{M}$ & turbine time constant & $\omega_0$ & sync. machine speed\\
\hline
$K_M$ & turbine gain & $\theta$ & rotor angle\\
\hline
$T_E$ & governor time constant & $t_{ij}$ & tie-line stiffness\\
\hline
$R$ & regulation constant & $P_E$ & steam valve position\\
\hline
$P_m$ & mechanical power & $K_E$ & speed governor gain\\
\hline
\end{tabular}
\end{table}
\normalsize

The agents aim to find the optimal mechanical power $P_m$ and rotor angle $\theta$ that minimize their operating cost, and at the same time, satisfy the coupled constraints on load balance. All these objectives are encoded by the following OPF problem which is adopted from Section 8.3.1 of \cite{JZ:2015}:
\begin{align}
\label{OPF formulation}
&\!\!\!\!\!\!\!\!\!\!\!\mathop {\min }\nolimits_{(P_m,\theta) \in {\mathbb{R}^{2N}}} \sum\nolimits_{i \in \mathcal{V}} {(a_iP_{m_i}^2+b_iP_{m_i}+c_i\theta_i^2)}\nonumber\\
{\rm{s.t.}}\,\,&\sum\nolimits_{j\in\mathcal{N}_i^{PH}}\frac{\omega_0t_{ij}}{2T_{H_i}}(\theta_{ij}-\sin\theta_{ij}(0))\nonumber\\
&=-\frac{D_i}{2T_{H_i}}\tilde\omega_i+\frac{\omega_0}{2T_{H_i}}(P_{m_i}-P_{m_i0}),\;\;\forall i\in\mathcal{V}.
\end{align}
In problem \eqref{OPF formulation}, the constraint uses the linearization $\sin\theta_{ij}\approx \theta_{ij}$ around the equilibrium.

In traditional hierarchical control, the decision-making layer first solves problem \eqref{OPF formulation} to obtain an optimal operating point $(P_m^*,\theta^*)$. After the convergence of the decision-making layer, the physical layer performs frequency control to regulate the physical dynamics \eqref{motivating example dynamics} to the desired operating point \cite{MH-DX-LW:2016}. To address the issues of hierarchical control, we aim to design an economic control algorithm such that the frequency control in the physical dynamics \eqref{motivating example dynamics} dynamically keeps track of the intermediate solutions of problem \eqref{OPF formulation}.


%
%


\section{Problem statement\label{section problem statement}}

This section provides the general formulation of the distributed economic control problem for dynamically coupled networks considered in this paper.

\subsection{Physical dynamics\label{problem formulation}}

Consider a set of agents $\mathcal{V} = \{ 1, \cdots ,N\}$.
The inherent physical dynamics is given by
\begin{align}
\label{physical dynamic}
\dot x_i (t) = \Xi_i(x(t),u_i(t))
\end{align}
where $x_i\in \mathbb{R}^{n_i}$ and $u_i\in\mathbb{R}$ are the state and control of agent $i$, respectively, with $n_i$ the dimension of $x_i$.
Let $x=(x_i)_{i=1}^N$, $u=(u_i)_{i=1}^N$ and $n=\sum\nolimits_{i=1}^N{n_i}$. 
The correspondence between \eqref{physical dynamic} and \eqref{thermal dynamics} (resp. \eqref{motivating example dynamics}) is given in Section \ref{OPF simulation} (resp. Section \ref{section simulation problem OPF}).

\subsection{Game formulation}

The control objective is formulated by a generalized Nash equilibrium problem. For each $i\in\mathcal{V}$, $w_i=[x_i^T,u_i]^T$ is agent $i$'s decision variable and is subject to a local constraint $W_i\subseteq\mathbb{R}^{r_i}$ with $r_i=n_i+1$. Let $w=(w_i)_{i=1}^N$, $r=\sum_{i=1}^Nr_i$ and $W=\prod_{i=1}^NW_i$.
Given $w_{-i}\in W_{-i}=\prod_{j\neq i}W_j$, each agent $i\in\mathcal{V}$ aims to solve the following optimization problem:
\begin{align}
\label{original optimization based game}
\min\nolimits_{w_i\in W_i}f_i(w)\quad{\rm s.t.}\;H(w)=0_{m},\;G(w)\leq0_{l}
\end{align}
where $f_i:W\to\mathbb{R}$ is the objective function of agent $i$, $H:W\to\mathbb{R}^m$ and $G:W\to\mathbb{R}^l$ are the functions of the coupled equality and inequality constraints shared by the agents. The correspondence between \eqref{original optimization based game} and \eqref{thermal OP} (resp. \eqref{OPF formulation}) is given in Section \ref{OPF simulation} (resp. Section \ref{section simulation problem OPF}).


Given any $w_{-i}\in W_{-i}$, the feasible set of $w_i$ is depicted by a set-valued map $\mathcal{D}_i: W_{-i}\to 2^{W_i}$ defined as ${\mathcal{D}_i}({w_{ - i}}) = \{ w_i\in W_i:H(w)=0_{m},\,G(w)\leq0_{l}\}$. With the map $\mathcal{D}_i$, problem \eqref{original optimization based game} can be equivalently expressed as follows:
\begin{align}
\label{general OP}
\mathop {\min }\nolimits_{w_i \in {\mathcal{D}_i}({w_{ - i}})} {f_i}(w_i,w_{ - i}).
\end{align}

Let $F=(f_i)_{i=1}^N$, $\mathcal{D}(w)=\prod_{i=1}^N\mathcal{D}_i(w_{-i})$ and $\mathcal{\hat D}=\{w\in W:H(w)=0_m,\,G(w)\leq0_l\}$. Notice that $\mathcal{D}(\cdot)$ is a set-valued map while $\mathcal{\hat D}$ is a constant set that represents the overall feasible set of $w$. The collection of \eqref{general OP} for all $i\in\mathcal{V}$ is a generalized Nash equilibrium problem with shared constraints \cite{Facchinei.Kanzow:07}, denoted by GNEP$(\mathcal{D},F)$. In particular, the players of the GNEP are the agents in $\mathcal{V}$, and the action profile of each player $i\in\mathcal{V}$ is $w_i\in\mathcal{D}_i(w_{-i})$. The solution concept of the GNEP is given by the following notion of Nash equilibrium (NE), at which no player can benefit from unilateral deviation.

\bdefinition
A joint state $\tilde w=(\tilde w_i)_{i=1}^N\in W$ is an NE of the GNEP$(\mathcal{D},F)$ if, for all $i\in\mathcal{V}$, $\tilde w_i \in \mathcal{D}_i(\tilde w_{-i})$ and ${f_i}(\tilde w) \leq {f_i}({w_i},\tilde w_{ - i})$ for any ${w_i}\in \mathcal{D}_i(\tilde w_{-i})$.
\label{def: constrained game with dynamics NE definition}
\edefinition

Let $\nabla F(w) = (\nabla_{w_i}f_i(w))_{i=1}^N$, $\nabla H(w)=(\nabla_{w}H_j(w))_{j=1}^m$ and $\nabla G(w)=(\nabla_{w}G_j(w))_{j=1}^l$. The following assumptions on the GNEP are supposed to hold throughout the paper.

\basm
For all $i\in\mathcal{V}$, for any $w_{-i}\in W_{-i}$, $f_i(\cdot,w_{-i})$ is of class $C^1$ and convex in $W_i$; $H(\cdot)$ is affine in $W$; $G(\cdot)$ is of class $C^1$ and convex in $W$.
\label{asm: convexity}
\easm

\basm
$\forall i\in\mathcal{V}$, $W_i$ is closed and convex.
\label{asm: convex feasible set}
\easm

\basm
The Slater condition holds for $\mathcal{\hat D}$, i.e., $\exists\bar w\in {\rm int}(W)$ such that $H(\bar w)=0_m$ and $G(\bar w)<0_l$.
\label{asm: Slater condition}
\easm

\basm
The maps $\nabla G$ and $\nabla F$ are locally Lipschitz on $W$.
\label{asm: Lipschitz continuity}
\easm

\basm
$\nabla F$ is strictly monotone over $W$.
\label{asm: strict monotone}
\easm

\bremark
\label{remark: game assumptions}
Assumptions \ref{asm: convexity}--\ref{asm: strict monotone} are standard in game analysis \cite{MZ-EF:Auto14,UR-UVS:2011,SL-AN:2013}. When $N=1$, i.e., there is a single agent, the GNEP reduces to an optimization problem. Assumptions \ref{asm: convexity}--\ref{asm: Lipschitz continuity} are standard to establish convergence for optimization problems. In particular, Assumptions \ref{asm: convexity} and \ref{asm: convex feasible set} lead to a convex optimization problem;
Assumption \ref{asm: Slater condition} implies that strong duality holds, i.e., the duality gap is zero; Assumption \ref{asm: Lipschitz continuity} implies that $\nabla G$ and $\nabla F$ are regular and is needed to establish the existence and uniqueness of the evolution of certain projected dynamics. When extending optimization problems to games, Assumption \ref{asm: strict monotone} is further needed to restrict the coupling between the agents.
\eremark

\subsection{From GNEP to variational inequality\label{section VE}}

A GNEP with shared constraints is in general difficult to solve, because it is defined over the joint set-valued map $\mathcal{D}=\{w\in\mathbb{R}^r:w_i\in\mathcal{D}_i(w_{-i}),\,\forall i\in\mathcal{V}\}$, which is in general non-convex in $w$. However, certain types of solutions can be calculated relatively easily by using the variational inequality (VI) approach \cite{GS-ea:10}. This subsection reviews relevant results of VI in the literature. We refer to \cite{FF-JSP:03,PE:10} for a comprehensive discussion on VI.

Given a set $\mathcal{Q}\subseteq \mathbb{R}^r$ and a map $\Phi:\mathcal{Q}\to\mathbb{R}^r$, the generalized variational inequality, denoted by GVI$(\mathcal{Q},\Phi)$, is the problem of finding $\tilde w\in\mathcal{Q}$ such that $(w-\tilde w)^T\Phi(\tilde w)\geq0$ for all $w\in\mathcal{Q}$. We next introduce another solution concept of the GNEP$(\mathcal{D},F)$, termed as variational equilibrium (VE).

\bdefinition
A joint state $\tilde w$ is a VE of the GNEP$(\mathcal{D},F)$ if it is a solution of the GVI$(\mathcal{\hat D},\nabla F)$.
\label{def: VE definition}
\edefinition

\bremark
If Assumptions \ref{asm: convexity} and \ref{asm: convex feasible set} hold, then any VE is an NE of the GNEP$(\mathcal{D},F)$, but not vice versa (\cite{PE:10}, Proposition 12.4).
There are two major motivations for considering VEs instead of NEs. First, the GVI is defined over the constant set $\hat{\mathcal{D}}$ which is convex in $w$. Hence, in general, the GVI is easier to solve than the GNEP. Second, in terms of optimality, VEs are in general more efficient than NEs. Specifically, an NE of a game is in general not optimal in terms of social welfare \cite{PD:1986}, except for some special cases, e.g., potential games \cite{DM-LSS:96}. 
Consider the following distributed separable optimization problem (DSOP): $\min_{w\in W}\sum\nolimits_{i =1}^N {{f_i}({w_i})}$, s.t. $H(w) = {0_m}$ and $G(w) \le {0_l}$. In general, an arbitrary NE of the GNEP$(\mathcal{D},F)$ with $\mathcal{D}$ and $F$ defined as above may not be an optimal solution of the DSOP. However, since the objective function is separable, the DSOP is equivalent to the GVI$(\mathcal{\hat D},\nabla F)$ (they share the same KKT conditions), and thus, by Definition \ref{def: VE definition}, the set of global optima of the DSOP is identical with the set of VEs of the GNEP$(\mathcal{D},F)$. The motivating problems in Section \ref{section motivation} are two instances of the DSOP.
\label{applicability to separable distributed optimization}
\eremark

In this paper, we focus on identifying VEs of the GNEP$(\mathcal{D},F)$, rather than general NEs. The issue of VE existence has been extensively studied. We refer to \cite{PE:10} for a thorough review of the existence issue. At this moment, we assume that the GNEP$(\mathcal{D},F)$ admits at least one VE. We will verify this assumption using suitable conditions when we consider specific GNEPs in the next two sections. With Assumptions \ref{asm: convexity}, \ref{asm: convex feasible set} and \ref{asm: strict monotone}, the GNEP$(\mathcal{D},F)$ has at most one VE (\cite{PE:10}, Proposition 12.9). Then given the existence assumption, the GNEP$(\mathcal{D},F)$ has a unique VE.



\subsection{Communication graph}

Denote by $\mathcal{G}=(\mathcal{V},\mathcal{E})$ the communication graph of the agents, where $\mathcal{E}\subseteq\mathcal{V}\times\mathcal{V}$ is the set of communication links such that $(i,j)\in\mathcal{E}$ if and only if agent $i$ can receive messages from agent $j$.
Denote by $\mathcal{N}_i\subseteq\mathcal{V}$ the set of neighbors of agent $i$ in $\mathcal{G}$, i.e., $\mathcal{N}_i=\{j\in\mathcal{V}\setminus\{i\}:(i,j)\in\mathcal{E}\}$.
%


\subsection{Objective\label{section objective}}

This paper aims to design distributed economic control algorithms such that $(x,u)$ of system \eqref{physical dynamic} keep track of the intermediate results of the GNEP$(\mathcal{D},F)$ all the time and globally asymptotically converge to the VE of the GNEP$(\mathcal{D},F)$.


\section{Case I: Linear dynamics\label{section linear dynamics}}

\subsection{Physical dynamics\label{subsection linear dynamics}}

In this section, we consider the case where system \eqref{physical dynamic} is in the following linear form:
\begin{align}
\label{linear physical dynamic}
\dot x_i (t) = A_{ii} x_i (t) + \sum\nolimits_{j \neq i} {{A_{ij}}{x_j}(t)} + B_i u_i (t).
\end{align}
If \eqref{linear physical dynamic} does not depend on some $x_j$, then $A_{ij}$ is a zero matrix. Let $A_i=[A_{i1},\cdots,A_{iN}]$, $A=(A_i)_{i=1}^N$ and $B=\mathrm{diag}(B_i)_{i=1}^N$. As stated in the last section, we aim to steer system \eqref{linear physical dynamic} to the VE of the GNEP$(\mathcal{D},F)$, a steady state but unknown \emph{a priori}. Thus stabilizability is insufficient as it requires to know the steady state in advance. Instead, this section imposes the following controllability assumption.

\basm
$\forall i\in\mathcal{V}$, $(A_{ii}, B_i)$ is controllable.
\label{asm: CCF}
\easm

\subsection{Communication graph\label{section linear communication graph}}

For each $i\in\mathcal{V}$, let ${\rm index}_H^i=\{\ell:H_\ell\;{\rm depends\;on\;}w_i\}$, ${\rm index}_G^i=\{\ell:G_\ell\;{\rm depends\;on\;}w_i\}$, $H^i=(H_\ell)_{\ell\in{\rm index}_H^i}$ and $G^i=(G_\ell)_{\ell\in{\rm index}_G^i}$.
That is, $H^i$ (resp. $G^i$) is the sub-vector of $H$ (resp. $G$) that depends on $w_i$. Denote the dimensions of $H^i$ and $G^i$ by $m_i$ and $l_i$, respectively. Notice that agent $i$ can drop those constraints $H_\ell(w)=0$ (resp. $G_\ell(w)\leq0$) such that $\ell\notin{\rm index}_H^i$ (resp. $\ell\notin{\rm index}_G^i$) from its optimization problem \eqref{general OP}, since these constraints do not affect $w_i$.
In this section, we have the following assumption on the communication graph $\mathcal{G}$, which states that $\mathcal{G}$ includes the dependency graph defined by \eqref{general OP} and \eqref{linear physical dynamic} as a subgraph. That is, agent $i$ needs to receive messages from agent $j$ whose $w_j$ affects the optimization problem \eqref{general OP} and/or the physical dynamics \eqref{linear physical dynamic} of agent $i$.

\basm
\label{asm: general communication graph}
For all $i\in\mathcal{V}$, if $[(A_ix)^T,f_i,H^{iT},G^{iT}]$ depends on $w_j$ for any $j\in\mathcal{V}\backslash\{i\}$, then $j\in\mathcal{N}_i$.
\easm

\subsection{Algorithm design\label{linear dynamic properties}}

In this section, we achieve the objective claimed in Section \ref{section objective} by an indirect approach. Roughly speaking: 1) we design an update rule for auxiliary variables such that they globally asymptotically converge to the VE, however, their evolution is not restricted by system \eqref{linear physical dynamic}; 2) for the physical variables $(x,u)$, we design a distributed control law for $u$ such that, governed by system \eqref{linear physical dynamic}, $(x,u)$ globally asymptotically keep track of the auxiliary variables and thus also asymptotically converge to the VE. In this section, we impose the following assumptions on the GNEP.

\basm
For all $i\in\mathcal{V}$, $W_i$ is bounded.
\label{asm: W bounded}
\easm

\basm
The equality constraint $Ax+Bu=0_n$ is included in $\hat{\mathcal{D}}$.
\label{asm: physical feasibility}
\easm

\bremark
Assumption \ref{asm: W bounded} establishes the existence of VE. Specifically, by Assumption \ref{asm: W bounded}, the set $\mathcal{\hat D}$ is bounded. Together with Assumptions \ref{asm: convexity} and \ref{asm: convex feasible set}, by Claim (8) of \cite{GS-ea:10}, the GNEP \eqref{general OP} admits a VE. Assumption \ref{asm: physical feasibility} guarantees that the VE of the GNEP satisfies the steady state condition of system \eqref{linear physical dynamic} so that the VE tracking problem is feasible. If $Ax+Bu=0_n$ is not originally included in $\hat{\mathcal{D}}$, one can deliberately add it into $\hat{\mathcal{D}}$. This assumption is also necessary for hierarchical control \cite{TH-XZ-WS-MZ-NL:2017CDC,BP-XZ-RS:2016}.
\eremark

\begin{algorithm2e}[htbp]
\caption{Distributed update rule and control law for linear system dynamics}
\label{Algorithm 1}

All the agents agree on any $\varepsilon>0$ s.t. $1/\varepsilon> 2\sqrt n (N - 1)\max_{i\in\mathcal{V}}\{\sigma_i \| {{P_{0i}}} \|_1\}+2$;

Each agent $i$ picks any $(\bar w_i(0),u_i(0))\in W_i\times\mathbb{R}$;

For any $\ell\in\{1,\cdots,m\}$ (resp. $\ell\in\{1,\cdots,l\}$), all agent $j$'s s.t. $\ell\in{\rm index}_H^j$ (resp. $\ell\in{\rm index}_G^j$) agree on any $\lambda_\ell(0)\in\mathbb{R}$ (resp. $\mu_\ell(0)\in\mathbb{R}_+$); each agent $i$ forms $\lambda^i(0)=(\lambda_\ell(0))_{\ell\in{\rm index}_H^i}$ and $\mu^i(0)=(\mu_\ell(0))_{\ell\in{\rm index}_G^i}$;
%

\While{$t\geq0$}{
Each agent $i$ sends $\bar w_i(t)$ to all $j$'s s.t. $i\in\mathcal{N}_j$;

Each agent $i$ updates its auxiliary variables by
\begin{align}
\label{decision update rule}
&\dot{\bar w}_i(t)=\Pi_{W_i}(\bar w_i(t),-\nabla_{\bar w_i}L_i(\bar w(t),\lambda^i(t),\mu^i(t))),\nonumber\\
&\dot\lambda^i(t)=H^i(\bar w(t)),\dot\mu^i(t)=\Pi_{\mathbb{R}_+^{l_i}}(\mu^i(t),G^i(\bar w(t))),
\end{align}

and executes the control law
\begin{align}
\label{control law for linear system}
u_i (t) = (K_i(\varepsilon)+a_i)\mathcal{T}_i(x_i(t)-\bar x_i(t))+\bar u_i(t).
\end{align}
}
\end{algorithm2e}
\normalsize

The algorithm is presented by Algorithm \ref{Algorithm 1} and detailed next. First we illustrate the design of the update rule \eqref{decision update rule} for auxiliary variables. For each $i\in\mathcal{V}$, define auxiliary variables $\bar x_i\in\mathbb{R}^{n_i}$ and $\bar u_i\in\mathbb{R}$, and let $\bar w_i=[\bar x_i^T,\bar u_i]^T$. Let $\bar x=(\bar x_i)_{i=1}^N$, $\bar u=(\bar u_i)_{i=1}^N$ and $\bar w=(\bar w_i)_{i=1}^N$. To handle constraints, for each $i\in\mathcal{V}$, define Lagrangian $L_i:W\times\mathbb{R}^{m_i}\times\mathbb{R}_+^{l_i}\to\mathbb{R}$ as $L_i({\bar w},\lambda^i,\mu^i) = {f_i}(\bar w)+\lambda^{iT}H^i(\bar w)+\mu^{iT}G^i(\bar w)$, where $\lambda^i\in\mathbb{R}^{m_i}$ and $\mu^i\in\mathbb{R}_+^{l_i}$ are agent $i$'s local copies of Lagrange multipliers $(\lambda_\ell)_{\ell\in{\rm index}_H^i}$ and $(\mu_\ell)_{\ell\in{\rm index}_G^i}$ associated with $H^i$ and $G^i$, respectively.
Each agent $i$ sends its real-time auxiliary variable $\bar w_i(t)$ to all agent $j$'s such that $i\in\mathcal{N}_j$.
This enables the agents to implement the update rule \eqref{decision update rule} for the auxiliary variables. For each $i\in\mathcal{V}$, the primal state $\bar w_i$ moves along the direction of $-\nabla_{\bar w_i}L_i(\bar w,\lambda^i,\mu^i)$ so as to minimize $L_i$, and the projection $\Pi_{W_i}$ guarantees $\bar w_i(t)\in W_i$ for all $t$. The dynamics of the dual states $(\lambda^i,\mu^i)$ are defined by the positive gradients of $L_i$ so as to maximize $L_i$. The projection $\Pi_{\mathbb{R}_+^{l_i}}$ on $\mu^i$ is to guarantee $\mu^i(t)\in\mathbb{R}_+^{l_i}$ for all $t$.

%

We next illustrate the design of the control law \eqref{control law for linear system} for the physical variables. The following preliminaries are needed to introduce the controller design. By Theorem 8.2 of \cite{chen}, for each $i\in\mathcal{V}$, since $(A_{ii},B_i)$ is controllable, there exists a nonsingular matrix $\mathcal{T}_i\in \mathbb{R}^{n_i\times n_i}$ such that $(\bar A_{ii},\hat B_i)$ with $\bar A_{ii}=\mathcal{T}_iA_{ii}\mathcal{T}_i^{-1}$ and $\hat B_i=\mathcal{T}_iB_i$ is in the controllable canonical form, i.e., $\bar A_{ii}=[[0_{n_i-1},I_{n_i-1}]^T,-a_i^T]^T$ with $a_i=[a_{i1},\cdots,a_{in_i}]$ a constant row vector, and $\hat B_i=[0_{n_i-1}^T,1]^T$.
Let $\hat x_i=\mathcal{T}_ix_i$, $v_i=u_i-a_i\hat x_i$, $\hat A_{ij}=\mathcal{T}_iA_{ij}\mathcal{T}_j^{-1}$ for all $j\neq i$, and $\hat A_{ii}=[[0_{n_i-1},I_{n_i-1}]^T,0_{n_i}]^T$. We then have $\dot{\hat x}_i(t)=\mathcal{T}_i\dot x_i(t)=\mathcal{T}_i(A_{ii}x_i(t)+\sum\nolimits_{j\neq i}A_{ij}x_j(t)+B_iu_i(t))=\bar A_{ii}\hat x_i(t)+\sum\nolimits_{j\neq i}\hat A_{ij}\hat x_j(t)+\hat B_iu_i(t)=\hat A_{ii}\hat x_i(t)+\sum\nolimits_{j\neq i}\hat A_{ij}\hat x_j(t)+\hat B_iv_i(t)$. For each $i\in\mathcal{V}$, notice that $(\hat A_{ii},\hat B_i)$ is controllable. Then there exists $K_{0i}\in \mathbb{R}^{1\times n_i}$ such that $\hat A_{ii}+\hat B_iK_{0i}$ is Hurwitz (\cite{chen}, Theorem 8.3). We further have that there exists a symmetric positive definite matrix $P_{0i}\in \mathbb{R}^{n_i\times n_i}$ such that $P_{0i}(\hat A_{ii}+\hat B_iK_{0i}) + (\hat A_{ii}+\hat B_iK_{0i})^TP_{0i} = -I_{n_i}$ (\cite{chen}, Theorem 5.5). For each $i\in\mathcal{V}$, let ${S_i}(\varepsilon ) \buildrel \Delta \over = \mathrm{diag}(\varepsilon^{j-1})_{j=1}^{n_i}$ and ${{\bar S}_i}(\varepsilon )\buildrel \Delta \over = \mathrm{diag}(\varepsilon^{j-n_i})_{j=1}^{n_i}$, where $\varepsilon\in \mathbb{R}$ is a constant scalar such that $0<\varepsilon<1$. The introduction of $\varepsilon$ is to enable arbitrarily fast tracking (by choosing $\varepsilon$ arbitrarily small). Let $S(\varepsilon)=\mathrm{diag}(S_i(\varepsilon))_{i=1}^N$, $K_i(\varepsilon) = \frac{1}{\varepsilon} K_{0i} \bar S_i(\varepsilon)$, $P_0=\mathrm{diag}(P_{0i})_{i=1}^N$ and $K(\varepsilon)=\mathrm{diag}(K_i(\varepsilon))_{i=1}^N$.

The distributed control law is given by \eqref{control law for linear system}. Each agent $i$ takes $(\bar x_i,\bar u_i)$ as an intermediate estimate of the VE and aims to drive $x_i$ by $u_i$ to asymptotically keep track of $\bar x_i$ and $u_i$ itself asymptotically keep track of $\bar u_i$.
Notice that \eqref{control law for linear system} does not require inter-agent communications.
We next informally explain how \eqref{control law for linear system} achieves asymptotic tracking. The formal proof is provided in Section \ref{appendix 2} using Laypunov analysis.

For each $i\in\mathcal{V}$, let $e_i=\mathcal{T}_i(x_i-\bar x_i)$ and $\hat e_i=S_i(\varepsilon)e_i$. To achieve asymptotic tracking, we aim that the error signal $\hat e_i$ converges to zero. By the controller \eqref{control law for linear system}, one can derive that the error dynamics is $\dot{\hat e}_i=\frac{1}{\varepsilon}(\hat A_{ii}+\hat B_iK_{0i})\hat e_i+S_i(\varepsilon)\sum_{j\neq i}\hat A_{ij}S_j^{-1}(\varepsilon)\hat e_j+S_i(\varepsilon)c_i$, where $c_i=\hat A_i\mathcal{T}_i\bar x_i+\hat B_i(\bar u_i-a_i\mathcal{T}_i\bar x_i)-\mathcal{T}_i\dot{\bar x}_i$ is free of $\varepsilon$. Recall that $\hat A_{ii}+\hat B_iK_{0i}$ is Hurwitz. Intuitively, by choosing arbitrarily small $\varepsilon>0$, the stability margin provided by the first term $\frac{1}{\varepsilon}(\hat A_{ii}+\hat B_iK_{0i})\hat e_i$ can be arbitrarily large, and the magnitude of the third term $S_i(\varepsilon)c_i$ can be arbitrarily small and thus dominated by the first term. In this section, we have the following assumption on the physical dynamics.


\basm
For all $i\in\mathcal{V}$, there exists a constant scalar $\sigma_i>0$ such that $\| S_i (\varepsilon )\hat A_{ij} S_j^{ - 1} (\varepsilon )\|_1 < \sigma_i$ for any $ j\in \mathcal{V}$ with $j\ne i$ and any $0<\varepsilon<1$.
\label{asm: S constraint}
\easm

\bremark
Notice that the second term of the error dynamics $\dot{\hat e}_i$, $S_i(\varepsilon)\sum_{j\neq i}\hat A_{ij}S_j^{-1}(\varepsilon)\hat e_j$, includes $S_j^{-1}(\varepsilon)$, which tends to infinity when $\varepsilon$ tends to zero. Assumption \ref{asm: S constraint} states that the coupling effects of $S_i(\varepsilon)\sum_{j\neq i}\hat A_{ij}S_j^{-1}(\varepsilon)$ is bounded for any $0<\varepsilon<1$ and is needed to achieve arbitrarily large stability margin. For any $i,j\in \mathcal{V}$, the entry $\bar a_{sc}$ of $S_i (\varepsilon )\hat A_{ij} S_j^{ - 1} (\varepsilon )$ at the $s$-th row and $c$-th column is $\varepsilon^{s-c}a_{sc}$, where $a_{sc}$ is the corresponding entry of $\hat A_{ij}$. A sufficient condition for Assumption \ref{asm: S constraint} is that $a_{sc}=0$ for all $s<c$. Under this condition, $\bar a_{sc}=0$ for all $s<c$. For the entries with $s\geq c$, since $0<\varepsilon<1$, $|\bar a_{sc}|\leq|a_{sc}|$. Thus, $\|S_i (\varepsilon )\hat A_{ij} S_j^{ - 1}$ $(\varepsilon )\|_1\leq \|\hat A_{ij}\|_1$ for any $0<\varepsilon<1$. If $n_i=n_j$, this condition means that $\hat A_{ij}$ is a lower triangular matrix. If $n_i>n_j$, $\hat A_{ij}$ is denser than lower triangular; if $n_i<n_j$, it is sparser.
\eremark

In lines 1 and 3 of Algorithm \ref{Algorithm 1}, the agreement of $\varepsilon$ and $\lambda(0)$ between the agents can be distributively carried out by the max-consensus scheme \cite{ROS-JAF-RMM:07}.

\subsection{Convergence results}

Theorem \ref{theorem: Algorithm 1} states that Algorithm \ref{Algorithm 1} drives the true states and inputs $(x(t),u(t))$ of system \eqref{linear physical dynamic} to globally asymptotically converge to the unique VE of problem \eqref{general OP}. The proof of Theorem \ref{theorem: Algorithm 1} is given in Section \ref{appendix 2}. 
In Section \ref{OPF simulation}, we will verify that the multi-zone building temperature regulation problem in Section \ref{section MZB} satisfies all the assumptions of Theorem \ref{theorem: Algorithm 1}.

\btheorem
Suppose Assumptions \ref{asm: convexity}--\ref{asm: strict monotone} and \ref{asm: CCF}--\ref{asm: S constraint} hold.
By Algorithm \ref{Algorithm 1}, $(x(t),u(t))$ of linear system \eqref{linear physical dynamic} globally asymptotically converge to the unique VE of problem \eqref{general OP}.
\label{theorem: Algorithm 1}
\etheorem

%
%
%
%

\section{Case II: Nonlinear dynamics\label{nonlinear section}}

\subsection{Physical dynamics\label{section nonlinear dynamics}}

Notice that Algorithm \ref{Algorithm 1} is inapplicable when system \eqref{physical dynamic} is nonlinear.
First, the steady state conditions of a nonlinear system introduce nonlinear functions into the equality constraints, making problem \eqref{general OP} non-convex. Second, the distributed controller \eqref{control law for linear system} cannot be applied to nonlinear systems.
If a linear system is stabilizable, there are many systematic ways to find its stabilizing controllers, e.g., cyclic design, Lyapunov-equation method and canonical-form method \cite{chen}. However, such general methods do not exist for nonlinear systems. In the literature of nonlinear control, researchers have studied stabilization of several classes of nonlinear systems with certain special structures. Notable examples include strict-feedback systems (or systems in lower triangular form) which can be stabilized by the technique of backstepping, and strict-feedforward systems (or systems in upper triangular form) which can be stabilized by the technique of nested saturation or forwarding \cite{AI:99,Khalil:02}. 
In this section,
%
we assume that system \eqref{physical dynamic} is nonlinear and in the following strict-feedback or lower-triangular form:
\begin{align}
\label{nonlinear state dynamics}
&\dot x_{i\ell}=\Gamma_{i\ell}(x^{[1\sim\ell]})+\Theta_{i\ell}(x^{[1\sim\ell]})x_{i(\ell+1)},\;\;\ell=1,\cdots,\bar n-1,\nonumber\\
&\dot x_{i\bar n}=\Gamma_{i\bar n}(x)+\Theta_{i\bar n}(x)u_i
\end{align}
where $x^{[1\sim\ell]}=[x^{[1]T},\cdots,x^{[\ell]T}]^T$ with $x^{[\ell]}=(x_{i\ell})_{i=1}^N$, and $\Gamma_{i\ell}$ and $\Theta_{i\ell}$ are in general nonlinear in $x^{[1\sim\ell]}$. System \eqref{nonlinear state dynamics} covers broad applications, e.g., multi-zone building systems \cite{ESK-MJY-CL:2008} and unmanned autonomous vehicles \cite{BA-HRP-MG:2010}. 


\subsection{Game reformulation\label{section game reformulation}}

Notice that the VE of the GNEP needs to satisfy the steady state condition of system \eqref{nonlinear state dynamics} so that the VE tracking problem is feasible. Since system \eqref{nonlinear state dynamics} is nonlinear, we cannot include its steady state condition into the equality constraints of GNEP. Otherwise, the affinity assumption on the equality constraints is violated. As illustrated next, this challenge can be overcome via a coordinate transformation which explicitly leverages the strict-feedback form of system \eqref{nonlinear state dynamics}.

For each $i\in\mathcal{V}$, the steady state condition of $x_{i1}$ is $\Gamma_{i1}(x^{[1]})+\Theta_{i1}(x^{[1]})x_{i2}=0$, by which we obtain the steady state of $x_{i2}$ as a function of $x^{[1]}$ as $\bar x_{i2}=-\Gamma_{i1}(x^{[1]})/\Theta_{i1}(x^{[1]})$.
The steady state condition of $x_{i2}$ is $\Gamma_{i2}(x^{[1\sim2]})+\Theta_{i2}(x^{[1\sim2]})x_{i3}=0$, by which we can obtain the steady state of $x_{i3}$ as $\bar x_{i3}=-\Gamma_{i2}(x^{[1\sim2]})/\Theta_{i2}(x^{[1\sim2]})$, which can be expressed as a function of $x^{[1]}$ by plugging in $\bar x_{j2}=-\Gamma_{j1}(x^{[1]})/\Theta_{j1}(x^{[1]})$ for all $j\in\mathcal{V}$.
Go on with this procedure of change of variables in the order of $x_{i2}\to\cdots\to x_{i\bar n}\to u_i$ to obtain the steady state of each $x_{i\ell}$ and $u_i$ as a function of only $x^{[1]}$, denoted by $\bar x_{i\ell}$ and $\bar u_i$, respectively. We then replace $(x_{i2},\cdots,x_{i\bar n},u_i)$ by $(\bar x_{i2},\cdots,\bar x_{i\bar n},\bar u_i)$ for all $i\in\mathcal{V}$ in the GNEP. The transformed GNEP then only depends on $x^{[1]}$. To slightly abuse the notations, we keep $f_i$, $W_i$, $H$ and $G$ as the functions associated with the transformed GNEP, which can be written as
\begin{align}
\label{OP nonlinear}
\mathop {\min }\nolimits_{x_{i1} \in {\mathcal{D}_i}({x_{ - i}^{[1]}})} {f_i}(x_{i1},x_{-i}^{[1]})
\end{align}
with ${\mathcal{D}_i}({x_{ - i}^{[1]}})=\{x_{i1}\in W_i:H^i(x_{i1},x_{-i}^{[1]})=0_{m_i},\,G^i(x_{i1},$ $x_{-i}^{[1]})\leq0_{l_i}\}$, where $H^i$ and $G^i$ are defined in Section \ref{section linear communication graph}.

For the transformed GNEP \eqref{OP nonlinear}, given the VE $\tilde x^{[1]}$, for each $i$, one can uniquely determine the steady states of $x_{i2},\cdots,x_{i\bar n},u_i$ by plugging $\tilde x^{[1]}$ into $\bar x_{i2},\cdots,\bar x_{i\bar n},\bar u_i$. Hence, it is guaranteed that the VE tracking problem is feasible.

\subsection{Algorithm design\label{section nonlinear algorithm design}}

Different from the indirect method of Algorithm \ref{Algorithm 1}, in this section, we directly work on the physical variables $(x,u)$ of system \eqref{nonlinear state dynamics} without introducing auxiliary variables, and design a controller $u$ that stabilizes system \eqref{nonlinear state dynamics} and meanwhile the steady state of $x^{[1]}$ is the VE of the GNEP \eqref{OP nonlinear}.
Inspired by backstepping \cite{KKK:95}, we introduce a coordinate transformation to covert system \eqref{nonlinear state dynamics} into a primal-dual dynamics perturbed by a chain of stable integrators. Similar to backstepping (Section 14.3 of \cite{Khalil:02}), the coordinate transformation requires the following assumption, which guarantees that all the derivatives in \eqref{coordinate transformation} and \eqref{control law for nonlinear dynamics} are well-defined and is necessary for the implementation of Algorithm \ref{Algorithm 2} introduced later in this subsection.


\basm
For all $i\in\mathcal{V}$, $W_i=\mathbb{R}$, $f_i$ is of class $C^{\bar n}$, $\Theta_{i\ell}$ and $\Gamma_{i\ell}$ are of class $C^\ell$, and $\Theta_{i\ell}(y)\neq0$ for any $y\in\mathbb{R}^{\ell N}$, for all $\ell=1,\cdots,\bar n$; there are no coupled inequality constraints.
\label{asm: smoothness}
\easm



In the rest of this section, to simplify notations, we omit the time index $t$ and the arguments of the component functions.

For each $i\in\mathcal{V}$, define Lagrangian $L_i(x^{[1]},\lambda^i)=f_i(x^{[1]})+\lambda^{iT}H^i(x^{[1]})$, where $\lambda^i\in\mathbb{R}^{m_i}$ is agent $i$'s local copy of Lagrange multiplier $(\lambda_\ell)_{\ell\in{\rm index}_H^i}$ associated with $H^i$. To illustrate the idea, let us first consider the special case where $\bar{n} = 3$. To achieve the VE, we aim that each $\dot x_{i1}$ converges to $-k_1\nabla_{x_{i1}}L_i$ for some $k_1>0$ so as to minimize $L_i$. However, $\dot x_{i1}$ is governed by \eqref{nonlinear state dynamics}. We then aim that the error between $\dot x_{i1}$ and $-k_1\nabla_{x_{i1}}L_i$ asymptotically diminishes. 
For each $i\in\mathcal{V}$, let $z_{i1}=x_{i1}$. We have $\dot z_{i1}=\Gamma_{i1}+\Theta_{i1}x_{i2}=-k_{1}\nabla_{z_{i1}}L_i+z_{i2}$, where $z_{i2}=k_{1}\nabla_{z_{i1}}L_i+\Gamma_{i1}+\Theta_{i1}x_{i2}$. Notice that $z_{i2}$ is the error between $\dot x_{i1}$ and $-k_1\nabla_{x_{i1}}L_i$. Hence, we aim that $z_{i2}$ asymptotically diminishes. We have $\dot z_{i2}=\frac{d}{dt}(k_{1}\nabla_{z_{i1}}L_i+\Gamma_{i1}+\Theta_{i1}x_{i2})=\frac{d}{dt}(k_{1}\nabla_{z_{i1}}L_i+\Gamma_{i1})+x_{i2}\frac{d}{dt}\Theta_{i1}+\Theta_{i1}(\Gamma_{i2}+\Theta_{i2}x_{i3})
=-k_{i2}z_{i2}+z_{i3}$, where $z_{i3}=k_{i2}z_{i2}+\frac{d}{dt}(k_{1}\nabla_{z_{i1}}L_i+\Gamma_{i1})+x_{i2}\frac{d}{dt}\Theta_{i1}+\Theta_{i1}(\Gamma_{i2}+\Theta_{i2}x_{i3})$ with constant $k_{i2}>0$. Intuitively, if $z_{i3}$ asymptotically diminishes, then $\dot z_{i2}$ asymptotically converges to $-k_{i2}z_{i2}$, which implies that $z_{i2}$ asymptotically diminishes. Hence, we further aim that $z_{i3}$ asymptotically diminishes. We can obtain $\dot z_{i3}=\frac{d}{dt}(k_{i2}z_{i2}+\frac{d}{dt}(k_{1}\nabla_{z_{i1}}L_i+\Gamma_{i1})+x_{i2}\frac{d}{dt}\Theta_{i1}+\Theta_{i1}(\Gamma_{i2}+\Theta_{i2}x_{i3}))= -k_{i2}^2z_{i2}+k_{i2}z_{i3}+\frac{d^2}{dt^2}(k_{1}\nabla_{z_{i1}}L_i+\Gamma_{i1})+\frac{d}{dt}(x_{i2}\frac{d}{dt}\Theta_{i1}+\Theta_{i1}\Gamma_{i2})+x_{i3}\frac{d}{dt}\Theta_{i1}\Theta_{i2}+ \Theta_{i1}\Theta_{i2}\Gamma_{i3}+\Theta_{i1}\Theta_{i2}\Theta_{i3}u_i$. To have $z_{i3}$ asymptotically diminishing, we aim that $\dot z_{i3}=-k_{i3}z_{i3}$ with constant $k_{i3}>0$. We then have $u_i=-\frac{1}{\Theta_{i1}\Theta_{i2}\Theta_{i3}}[-k_{i2}^2z_{i2}+(k_{i2}+k_{i3})z_{i3}+\frac{d^2}{dt^2}(k_{1}\nabla_{z_{i1}}L_i+\Gamma_{i1})+\frac{d}{dt}(x_{i2}\frac{d}{dt}\Theta_{i1}+\Theta_{i1}\Gamma_{i2})+x_{i3}\frac{d}{dt}\Theta_{i1}\Theta_{i2}+ \Theta_{i1}\Theta_{i2}\Gamma_{i3}]$.

Following the above procedure, for the general case, sequentially define a coordinate transformation for $z_{i2},\cdots,z_{i\bar n}$ and design a controller $u_i$ such that system \eqref{nonlinear state dynamics} is converted into the following:
\begin{align}
\label{transformed states dynamics}
\dot z_{i1}&=-k_{1}\nabla_{z_{i1}}L_i+z_{i2},\nonumber\\
\dot z_{i\ell}&=-k_{i\ell}z_{i\ell}+z_{i(\ell+1)},\;\;\forall\ell=2,\cdots,\bar n-1,\nonumber\\
\dot z_{i\bar n}&=-k_{i\bar n}z_{i\bar n}
\end{align}
where $k_1,k_{i\ell},\cdots,k_{i\bar n}>0$. In this section, we have the following assumption.

\basm
The map $\nabla F$ is strongly monotone on $\mathbb{R}^N$ with constant $M>0$.
\label{asm: strong monotone}
\easm

\bremark
In system \eqref{transformed states dynamics}, the dynamics of $z_{i2},\cdots,z_{i\bar n}$ is a sequence of stable, cascading and single integrators. Thus, $z_{i2},\cdots,z_{i\bar n}$ are diminishing. Notice that the dynamics of $z_{i1}$ is the primal-dual dynamics perturbed by $z_{i2}$. Strict monotonicity of $\nabla F$ is not sufficient to ensure convergence because it provides zero stability margin for the primal-dual dynamics. Instead, the stronger property of strong monotonicity is needed to ensure nonzero stability margin, and then the convergence of Algorithm \ref{Algorithm 2}. Meanwhile, a byproduct of Assumption \ref{asm: strong monotone} is the existence and uniqueness of VE. Specifically, with Assumptions \ref{asm: convexity}, \ref{asm: convex feasible set} and \ref{asm: strong monotone}, by Theorem 2.3.3 of \cite{FF-JSP:03}, the GNEP \eqref{OP nonlinear} has a unique VE.
\eremark

We have illustrated the idea of controller design for system \eqref{nonlinear state dynamics} via a backstepping-inspired coordinate transformation. We next provide their closed forms, which are derived from \eqref{transformed states dynamics} by expanding $\dot z_{i1},\cdots,\dot z_{i\bar n}$ on its left-hand side.
For convenience of demonstration, we define the following notations. For any $\alpha,\beta\in\{1,\cdots,\bar n-1\}$, let $T_{\alpha,\beta}^i$ be the sum of all the $\alpha$-th order products of the first $\beta$ elements of $\{k_{i2},\cdots,k_{i\bar n}\}$ with repeat. For example, $T_{1,2}^i=k_{i2}+k_{i3}$ and $T_{2,3}=k_{i2}^2+k_{i3}^2+k_{i4}^2+k_{i2}k_{i3}+k_{i2}k_{i4}+k_{i3}k_{i4}$. Let $T_{\alpha,\beta}^i=0$ for any $\alpha\leq0$. For each $i\in\mathcal{V}$, for each $\ell\in\{1,\cdots,\bar n\}$, let $\Theta_{i}^{[\ell]} =\Theta_{i1}\Theta_{i2}\cdots\Theta_{i\ell}$. For each $i\in\mathcal{V}$, agent $i$ performs the following coordinate transformation:
\begin{align}
\label{coordinate transformation}
&z_{i1}=x_{i1},\;\;z_{i2}=k_{1}\nabla_{z_{i1}}L_i+\Theta_{i1}x_{i2}+\Gamma_{i1},\nonumber\\
&z_{i\ell}=(-1)^{\ell-1}T_{\ell-2,1}^iz_{i2}+(-1)^{\ell-2}T_{\ell-3,2}^iz_{i3}+\cdots\nonumber\\
&+(-1)^2T_{1,\ell-2}^iz_{i(\ell-1)}+\frac{d^{\ell-2}}{dt^{\ell-2}}(k_1\nabla_{z_{i1}}L_i+\Gamma_{i1})\nonumber\\
&+\frac{d^{\ell-3}}{dt^{\ell-3}}(x_{i2}\frac{d}{dt}\Theta_{i1}+\Theta_{i1}\Gamma_{i2})+\cdots+\frac{d}{dt}(x_{i(\ell-2)}\frac{d}{dt}\Theta_{i}^{[\ell-3]}\nonumber\\
&+\Theta_{i}^{[\ell-3]}\Gamma_{i(\ell-2)})+x_{i(\ell-1)}\frac{d}{dt}\Theta_{i}^{[\ell-2]}+\Theta_{i}^{[\ell-2]}\Gamma_{i(\ell-1)}\nonumber\\
&+\Theta_{i}^{[\ell-1]}x_{i\ell},\;\;\forall\ell=3,4,\cdots,\bar n.
\end{align}

For each $i\in\mathcal{V}$, agent $i$'s control $u_i$ is designed as follows:
\begin{align}
\label{control law for nonlinear dynamics}
&u_i=-\frac{1}{\Theta_{i}^{[\bar n]}}[(-1)^{\bar n}T_{\bar n-1,1}^iz_{i2}+(-1)^{\bar n-1}T_{\bar n-2,2}^iz_{i3}+\cdots\nonumber\\
&+(-1)^2T_{1,\bar n-1}^iz_{i\bar n}+\frac{d^{\bar n-1}}{dt^{\bar n-1}}(k_1\nabla_{z_{i1}}L_i+\Gamma_{i1})+\frac{d^{\bar n-2}}{dt^{\bar n-2}}\nonumber\\
&(x_{i2}\frac{d}{dt}\Theta_{i1}+\Theta_{i1}\Gamma_{i2})+\cdots+\frac{d}{dt}(x_{i(\bar n-1)}\frac{d}{dt}\Theta_{i}^{[\bar n-2]}\nonumber\\
&+\Theta_{i}^{[\bar n-2]}\Gamma_{i(\bar n-1)})+x_{i\bar n}\frac{d}{dt}\Theta_{i}^{[\bar n-1]}+\Theta_{i}^{[\bar n-1]}\Gamma_{i\bar n}].
\end{align}

With the strict-feedback form of \eqref{nonlinear state dynamics}, for each $i\in\mathcal{V}$ and each $\ell=1,\cdots,\bar n$, $z_{i\ell}$ only depends on $z^{[1\sim(\ell-1)]}$ but does not depend on $z^{[\ell']}$ for any $\ell'\geq\ell$, and $u_i$ only depends on $z$ but does not depend on any other $u_j$. Hence, \eqref{coordinate transformation} and \eqref{control law for nonlinear dynamics} can be sequentially computed in the order of $z_{i1}\to z_{i2}\to\cdots\to z_{i\bar n}\to u_i$. The coordinate transformation \eqref{coordinate transformation} and controller \eqref{control law for nonlinear dynamics} impose the following assumption on the communication graph $\mathcal{G}$, which guarantees that $\mathcal{G}$ includes the dependency graph defined by \eqref{coordinate transformation} and \eqref{control law for nonlinear dynamics} as a subgraph and is necessary for the implementation of Algorithm \ref{Algorithm 2}.

\basm
\label{asm: nonlinear communication graph}
For all $i\in\mathcal{V}$, if there exists $q\in\{0,1,\cdots,\bar n-1\}$ such that $\frac{d^q}{dt^q}[\Gamma_i^T,\Theta_i^T,f_i,H^{iT}]$ depends on $x_j$ for any $j\in\mathcal{V}\backslash\{i\}$, then $j\in\mathcal{N}_i$.
\easm


The following lemma formally verifies the dynamics of the transformed system derived by \eqref{coordinate transformation} and \eqref{control law for nonlinear dynamics}.

\blemma
Under Assumption \ref{asm: smoothness}, by \eqref{coordinate transformation} and \eqref{control law for nonlinear dynamics}, the dynamics of the transformed system follows \eqref{transformed states dynamics}.
\label{lemma: coordinate transformation dynamics}
\elemma


%

The algorithm is summarized by Algorithm \ref{Algorithm 2}. Each agent $i$ executes the control law \eqref{control law for nonlinear dynamics} and updates the transformed state $z_i$ by \eqref{transformed states dynamics}. The update rule for the dual variable $\lambda^i$ is similar to that in Algorithm \ref{Algorithm 1}.
In lines 1 and 5 of Algorithm \ref{Algorithm 2}, the agreement of $k_1$ and $\lambda(0)$ between the agents can be distributively carried out by the max-consensus scheme \cite{ROS-JAF-RMM:07}.

\begin{algorithm2e}[htbp]
\caption{Distributed control law for nonlinear system dynamics}
\label{Algorithm 2}

All the agents agree on any $k_1>0$ s.t. $k_1M>1$;

Each agent $i$ picks any $(x_i(0),u_i(0))\in \mathbb{R}^{\bar n}\times\mathbb{R}$, and any $k_{i\ell}>1$ for all $\ell=2,\cdots,\bar n$;

Each agent $i$ sends $x_i(0)$ to all $j$'s s.t. $i\in\mathcal{N}_j$;

Each agent $i$ computes $z_i(0)$ by \eqref{coordinate transformation};

For any $\ell\in\{1,\cdots,m\}$, all agent $j$'s s.t. $\ell\in{\rm index}_H^j$ agree on any $\lambda_\ell(0)\in\mathbb{R}$; each agent $i$ forms $\lambda^i(0)=(\lambda_\ell(0))_{\ell\in{\rm index}_H^i}$;
%

\While{$t\geq0$}{
Each agent $i$ sends $x_i(t)$ to all $j$'s s.t. $i\in\mathcal{N}_j$;

Each agent $i$ executes the control law \eqref{control law for nonlinear dynamics}, updates $z_i$ by \eqref{transformed states dynamics}, and updates $\lambda^i$ by
\begin{align}
\label{nonlinear dual update}
&\dot\lambda^i=k_1H^i(z^{[1]}).
\end{align}
}
\end{algorithm2e}

\normalsize

\subsection{Convergence results}

Theorem \ref{theorem: Algorithm 2} states that Algorithm \ref{Algorithm 2} drives the states and inputs $(x(t),u(t))$ to globally asymptotically converge to a steady state of system \eqref{nonlinear state dynamics} such that the limiting steady state of $x^{[1]}$ is the unique VE of problem \eqref{OP nonlinear}. The proof of Theorem \ref{theorem: Algorithm 2} is given in Section \ref{proof of theorem 2}. 
In Section \ref{section simulation problem OPF}, we will verify that the OPF problem in Section \ref{section OPF} satisfies all the assumptions of Theorem \ref{theorem: Algorithm 2}.

\btheorem
Suppose Assumptions \ref{asm: convexity}--\ref{asm: Lipschitz continuity} and \ref{asm: smoothness}--\ref{asm: nonlinear communication graph} hold.
By Algorithm \ref{Algorithm 2}, $(x(t),u(t))$ globally asymptotically converge to a steady state of system \eqref{nonlinear state dynamics}, denoted by $(\tilde x,\tilde u)$, with $\tilde x^{[1]}$ being the unique VE of problem \eqref{OP nonlinear}.
\label{theorem: Algorithm 2}
\etheorem


\subsection{Discussion}

The two problems in Section \ref{section linear dynamics} and Section \ref{nonlinear section} are complementary to each other, rather than one is a special case of the other. In Section \ref{section nonlinear dynamics}, we have illustrated that Algorithm \ref{Algorithm 1} may not be applicable to the problem in Section \ref{nonlinear section}. Here, we further clarify that Algorithm \ref{Algorithm 2} may not be applicable to the problem in Section \ref{section linear dynamics}. First, in Section \ref{nonlinear section}, it requires that the physical dynamics \eqref{nonlinear state dynamics} is in the strict-feedback form. 
However, the physical dynamics \eqref{linear physical dynamic} in Section \ref{section linear dynamics} is not necessarily strict-feedback. 
Second, for the GNEP of Section \ref{nonlinear section}, it requires that $W_i=\mathbb{R}$ and $f_i$ is of class $C^{\bar n}$ (please refer to Assumption \ref{asm: smoothness}). This assumption may not hold for the GNEP of Section \ref{section linear dynamics}, where $W_i$ is only assumed to be closed and convex (please refer to Assumption \ref{asm: convex feasible set}) and $f_i$ is only assumed to be of class $C^1$ (please refer to Assumption \ref{asm: convexity}). Third, for the GNEP of Section \ref{nonlinear section}, it requires that $\nabla F$ is strongly monotone (please refer to Assumption \ref{asm: strong monotone}), while for the GNEP of Section \ref{section linear dynamics}, it only requires that $\nabla F$ is strictly monotone (please refer to Assumption \ref{asm: strict monotone}).

A major difference between Algorithm \ref{Algorithm 1} and Algorithm \ref{Algorithm 2} is that the former needs auxiliary variables for the decision-making process, while the latter does not. The reason why Algorithm \ref{Algorithm 1} needs auxiliary variables is that the decision making dynamics \eqref{decision update rule} cannot be included in the physical dynamics \eqref{linear physical dynamic}. Hence, the intermediate estimates of the decision making dynamics \eqref{decision update rule} may not be equilibria of the physical dynamics \eqref{linear physical dynamic}. The auxiliary variables are introduced such that \eqref{linear physical dynamic} and \eqref{decision update rule} are cascaded. 
For the problem of Section \ref{nonlinear section}, thanks to the strict-feedback form of the physical dynamics \eqref{nonlinear state dynamics}, one can perform a coordinate transformation using the steady state conditions of \eqref{nonlinear state dynamics} such that the GNEP only depends on $x^{[1]}$ (please refer to Section \ref{section game reformulation}). Once $x^{[1]}$ is derived by the decision-making process, one can then sequentially determine the rest states and control inputs by reversing the coordinate transformation, which exactly follows the physical dynamics \eqref{nonlinear state dynamics}. In other words, for the nonlinear system case in Section \ref{nonlinear section}, the decision making dynamics is embedded in the physical dynamics, and the intermediate decision making estimates are always feasible to the physical dynamics. 
Hence, Algorithm \ref{Algorithm 2} does not need auxiliary variables.

\section{Proofs}

This section proves the results of Sections \ref{section linear dynamics} and \ref{nonlinear section}.

\subsection{Proof of Theorem \ref{theorem: Algorithm 1}\label{appendix 2}}

This subsection is devoted to proving Theorem \ref{theorem: Algorithm 1}.

\begin{IEEEproof} We first prove the convergence of the auxiliary variables $(\bar w,\lambda,\mu)$ in the following theorem.

\btheorem
\label{theorem auxiliary variable convergence}
Any solution of the primal-dual dynamics \eqref{decision update rule} starting from any point $(w(0),\lambda(0),\mu(0))\in W\times\mathbb{R}^m\times\mathbb{R}_+^l$ asymptotically converges to a point $(\tilde w,\tilde\lambda,\tilde\mu)$ where $\tilde w$ is the VE of problem \eqref{general OP}.
\etheorem

\begin{IEEEproof}
Let $\nabla_{\bar w}L(\bar w,\lambda,\mu)=(\nabla_{\bar w_i}L_i(\bar w,\lambda,\mu))_{i=1}^N$. Under Assumptions \ref{asm: convexity}, \ref{asm: convex feasible set} and \ref{asm: Slater condition}, by Section 12.2.3 of \cite{PE:10} and problem 5 on page 266 of \cite{DGL:1969}, $\tilde w\in W$ is the VE of the GNEP$(\mathcal{D},F)$ if and only if there exist $(\tilde \lambda,\tilde\mu)\in\mathbb{R}^m\times\mathbb{R}_+^l$ such that the KKT conditions hold
\begin{align}
\label{game KKT condition}
&\nabla_{\bar w}L(\tilde w,\tilde\lambda,\tilde\mu)\in N_W(\tilde w),\,\,H(\tilde w)=0_m,\nonumber\\
&G(\tilde w)\leq0_l\,\,{\rm and}\,\,\tilde \mu^{T}G(\tilde w)=0.
\end{align}
We note that the above KKT conditions only characterize VEs, rather than arbitrary NEs. For differences between the KKT conditions for VE and NE, please refer to Section 3 of \cite{FF-AF-VP:2007}.

A triple $(\tilde w,\tilde\lambda,\tilde\mu)$ satisfying \eqref{game KKT condition} is called a primal-dual optimizer of \eqref{general OP}. Denote by $Q$ the set of primal-dual optimizers of \eqref{general OP}. The update rule \eqref{decision update rule} can be written compactly as
\begin{align}
\label{update rule standard form}
&(\dot{\bar w},\dot\lambda,\dot\mu)=\Pi_{W\times\mathbb{R}^m\times\mathbb{R}_+^l}((\bar w,\lambda,\mu),\nonumber\\
&(-\nabla_{\bar w}L(\bar w,\lambda,\mu),H(\bar w),G(\bar w))).
\end{align}

Define the Lyapunov function as $V(\bar w,\lambda,\mu)=\frac{1}{2}(\|\bar w-\tilde w\|^2+\|\lambda-\tilde \lambda\|^2+\|\mu-\tilde \mu\|^2)$, where $(\tilde w,\tilde\lambda,\tilde\mu)$ is any point of $Q$. For convenience of notation, denote by $\mathcal{L}V(\bar w,\lambda,\mu)$ the Lie derivative of $V$ along \eqref{update rule standard form} at $(\bar w,\lambda,\mu)\in W\times\mathbb{R}^m\times\mathbb{R}_+^l$. We next show that system \eqref{update rule standard form} satisfies all the hypotheses of Lemma \ref{lemma invariance} by proving the following series of claims. First, the following claim establishes the result that $V$ is non-increasing.

\bclaim
\label{middle lemma monotonicity of Lyapunov function}
It holds that $\mathcal{L}V(\bar w,\lambda,\mu)\leq0$ for all $(\bar w,\lambda,\mu)\in W\times\mathbb{R}^m\times\mathbb{R}_+^l$.
\eclaim

\begin{IEEEproof}
We first show for all $(\bar w,\lambda,\mu)\in W\times\mathbb{R}^m\times\mathbb{R}_+^l$, $(\bar w-\tilde w)^T\Pi_W(\bar w,-\nabla_{\bar w}L(\bar w,\lambda,\mu))\leq-(\bar w-\tilde w)^T\nabla_{\bar w}L(\bar w,\lambda,\mu)$. Let $y=-\nabla_{\bar w}L(\bar w,\lambda,\mu)$. By Lemma \ref{lemma vector projection}, we only have to consider the case where $\bar w\in{\rm bd}(W)$ and $y^T\gamma^*(\bar w)<0$ with $\gamma^*(\bar w)={\rm argmin}_{\gamma\in N_W(\bar w)}y^T\gamma$. By Lemma \ref{lemma vector projection}, we have $\Pi_W(\bar w,y)=y-(y^T\gamma^*(\bar w))\gamma^*(\bar w)$. Then we have $(\bar w-\tilde w)^T\Pi_W(\bar w,y)=(\bar w-\tilde w)^T(y-(y^T\gamma^*(\bar w))\gamma^*(\bar w))=(\bar w-\tilde w)^Ty-(y^T\gamma^*(\bar w))(\bar w-\tilde w)^T\gamma^*(\bar w)$. Since $\bar w,\tilde w\in W$ and $\gamma^*(\bar w)\in N_W(\bar w)$, by the definition of normal cone $N_W$, we have $(\bar w-\tilde w)^T\gamma^*(\bar w)\leq0$. Since $y^T\gamma^*(\bar w)<0$, we have $(y^T\gamma^*(\bar w))(\bar w-\tilde w)^T\gamma^*(\bar w)\geq0$ and thus $(\bar w-\tilde w)^T\Pi_W(\bar w,y)\leq(\bar w-\tilde w)^Ty$. Following the similar argument, we can obtain $(\mu-\tilde \mu)^T\Pi_{\mathbb{R}_+^l}(\mu,G(\bar w))\leq(\mu-\tilde\mu)^TG(\bar w)$. By the definition of $N_W$ and the first equation of \eqref{game KKT condition}, we have $(\tilde w-\bar w)^T\nabla_{\bar w}L(\tilde w,\tilde\lambda,\tilde\mu)\leq0$ for any $\bar w\in W$.
By the affinity of $H$ and the convexity of $G$ (see Assumption \ref{asm: convexity}), we have
\begin{align}
\label{middle 1}
&-(\bar w-\tilde w)^T\nabla_{\bar w}L(\bar w,\lambda,\mu)\leq-(\bar w-\tilde w)^T(\nabla F(\bar w)\nonumber\\
&+\sum\nolimits_{j=1}^m\lambda_j\nabla H_j(\bar w)+\sum\nolimits_{j=1}^l\mu_j\nabla G_j(\bar w)-\nabla F(\tilde w)\nonumber\\
&-\sum\nolimits_{j=1}^m\tilde\lambda_j\nabla H_j(\tilde w)-\sum\nolimits_{j=1}^l\tilde\mu_j\nabla G_j(\tilde w))\nonumber\\
&\leq-(\bar w-\tilde w)^T(\nabla F(\bar w)-\nabla F(\tilde w))-(\lambda-\tilde\lambda)^T\nonumber\\
&\quad(H(\bar w)-H(\tilde w))-(\mu-\tilde\mu)^T(G(\bar w)-G(\tilde w)).
\end{align}

By \eqref{game KKT condition}, $H(\tilde w)=0_m$, $G(\tilde w)\leq0_l$ and $\tilde\mu^TG(\tilde w)=0$. Since $\mu\in\mathbb{R}_+^l$, $\mu^TG(\tilde w)\leq0$. By \eqref{middle 1}, $-(\bar w-\tilde w)^T\nabla_{\bar w}L(\bar w,\lambda,\mu)\leq-(\bar w-\tilde w)^T(\nabla F(\bar w)-\nabla F(\tilde w))-(\lambda-\tilde\lambda)^TH(\bar w)-(\mu-\tilde\mu)^TG(\bar w)$.
Then, by strict monotonicity of $\nabla F$, we have
\begin{align}
\label{middle 2}
&\mathcal{L}V(\bar w,\lambda,\mu)=(\bar w-\tilde w)^T\Pi_W(\bar w,-\nabla_{\bar w}L(\bar w,\lambda,\mu))\nonumber\\
&\quad+(\lambda-\tilde\lambda)^TH(\bar w)+(\mu-\tilde\mu)^T\Pi_{\mathbb{R}_+^l}(\mu,G(\bar w))\nonumber\\
&\leq-(\bar w-\tilde w)^T\nabla_{\bar w}L(\bar w,\lambda,\mu)+(\lambda-\tilde\lambda)^TH(\bar w)+(\mu-\tilde\mu)^T\nonumber\\
&\quad G(\bar w)\leq-(\bar w-\tilde w)^T(\nabla F(\bar w)-\nabla F(\tilde w))\leq0.
\end{align}

By strict monotonicity of $\nabla F$, the equality sign of the last equation of \eqref{middle 2} holds if and only if $\bar w=\tilde w$. Since $(\bar w,\lambda,\mu)$ is an arbitrary triple in $W\times\mathbb{R}^m\times\mathbb{R}_+^l$, we have $\mathcal{L}V(\bar w,\lambda,\mu)\leq0$ for all $(\bar w,\lambda,\mu)\in W\times\mathbb{R}^m\times\mathbb{R}_+^l$.
\end{IEEEproof}


The following claim establishes the existence and uniqueness of the solution of \eqref{decision update rule}. Moreover, the solution remains in a bounded set.

\bclaim
\label{middle lemma existence and uniqueness}
For any starting point $(\bar w(0),\lambda(0),\mu(0))\in W\times\mathbb{R}^m\times\mathbb{R}_+^l$, a unique solution $\Gamma:[0,\infty)\to W\times\mathbb{R}^m\times\mathbb{R}_+^l$ of \eqref{decision update rule} exists and $\Gamma(t)\in (W\times\mathbb{R}^m\times\mathbb{R}_+^l)\cap V^{-1}(\leq V(\bar w(0),\lambda(0),\mu(0)))$ for any $t\geq0$.
\eclaim

\begin{IEEEproof}
By Assumption \ref{asm: convexity}, $\nabla F$ and $\nabla G$ are all locally Lipschitz on $W$. Since $H$ is affine in $w$, $\nabla H$ is Lipschitz on $W$. The proof of Claim \ref{middle lemma existence and uniqueness} then follows the same argument of the proof of Lemma 4.3 of \cite{AC-EM-JC:2016}.
\end{IEEEproof}

The following claim establishes the invariance of the omega-limit set of any solution of \eqref{update rule standard form}.

\bclaim
\label{lemma invariance omega limit set}
The omega-limit set of any solution of \eqref{update rule standard form} is invariant under \eqref{update rule standard form}.
\eclaim

\begin{IEEEproof}
The proof of Claim \ref{lemma invariance omega limit set} follows the same argument of the proof of Lemma 4.1 of \cite{Khalil:02}.
\end{IEEEproof}

For any $\delta>0$, consider the compact set $S=V^{-1}(\leq \delta)\cap(W\times\mathbb{R}^m\times\mathbb{R}_+^l)$. By Claim \ref{middle lemma existence and uniqueness}, starting from any point in $S$, there exists a unique solution of \eqref{update rule standard form} and $S$ is invariant under \eqref{update rule standard form}. Furthermore, by Claim \ref{lemma invariance omega limit set}, the omega-limit set of each solution starting from any point in $S$ is invariant. Finally, by Claim \ref{middle lemma monotonicity of Lyapunov function}, $\mathcal{L}V(\bar w,\lambda,\mu)\leq0$ for all $(\bar w,\lambda,\mu)\in S$. By Lemma \ref{lemma invariance}, starting from any point in $S$, the solution of \eqref{update rule standard form} converges to the largest invariant set $M$ contained in ${\rm cl}(Z)$, where $Z=\{(\bar w,\lambda,\mu)\in S:\mathcal{L}V(\bar w,\lambda,\mu)=0\}$. By the proof of Claim \ref{middle lemma monotonicity of Lyapunov function}, $\mathcal{L}V(\bar w,\lambda,\mu)=0$ if and only if $\bar w=\tilde w$ and $\mu^TG(\tilde w)=0$. Hence, we have $Z=\{(\bar w,\lambda,\mu)\in S:\bar w=\tilde w,\;\mu^TG(\tilde w)=0\}$. It is clear that $Z$ is closed. Let $(\tilde w,\lambda,\mu)\in M\subseteq Z$. The solution of \eqref{update rule standard form} starting from $(\tilde w,\lambda,\mu)$ remains in $M$ only if $\nabla_{\bar w}L(\tilde w,\lambda,\mu)\in N_W(\tilde w)$. Together with $H(\tilde w)=0_m$, $G(\tilde w)\leq 0_l$ and $\mu^TG(\tilde w)=0$, we have that the triple $(\tilde w,\lambda,\mu)$ satisfies the KKT conditions \eqref{game KKT condition}. Hence, we have $M\subseteq Q$. Since $\delta$ is arbitrary, we have that $Q$ is globally asymptotically stable on $W\times\mathbb{R}^m\times\mathbb{R}_+^l$. Finally, by the definition of omega-limit set and Claim \ref{middle lemma monotonicity of Lyapunov function}, the omega-limit set of any solution of \eqref{update rule standard form} is a singleton. This implies that any solution of \eqref{update rule standard form} converges to a point in $Q$. This completes the proof of Theorem \ref{theorem auxiliary variable convergence}.
\end{IEEEproof}

Next we prove the convergence of the physical variables $(x,u)$. For each $i\in\mathcal{V}$, let $\bar{\bar x}_i=\mathcal{T}_i\bar x_i$, $\bar v_i=\bar u_i-a_i\bar{\bar x}_i$, $e_i=\hat x_i-\bar{\bar x}_i=\mathcal{T}_i(x_i-\bar x_i)$, $\hat e_i=S_i(\varepsilon)e_i$ and $e_i^v=v_i-\bar v_i=K_i(\varepsilon)e_i$. By \eqref{control law for linear system}, we have $e_i^v=K_i(\varepsilon)e_i$. Let $(e,\hat e,e^v)=(e_i,\hat e_i,e_i^v)_{i=1}^N$. By the definition of $S_i(\varepsilon)$ and $\bar S_i(\varepsilon)$, we can obtain $\varepsilon {S_i(\varepsilon)}{\hat A_{ii}} = {\hat A_{ii}}{S_i(\varepsilon)}$, ${S_i(\varepsilon)}{\hat B_i} = {\varepsilon ^{{n_i} - 1}}{\hat B_i}$ and ${\varepsilon ^{{n_i} - 1}}{{\bar S}_i(\varepsilon)} = {S_i(\varepsilon)}$. Choose the Lyapunov function as ${V}(\hat e) = {\hat e^T}P_0\hat e$. We then have
$\dot V(\hat e) = 2\hat e^T P_0\dot {\hat {e}}=2\hat e^TP_0S(\varepsilon)\dot e= 2\hat e^TP_0S(\varepsilon)(\hat Ae+\hat Be^v+\hat A{\bar {\bar x}}+\hat B\bar v-\dot {\bar{\bar x}})$. Let $c=\hat A{\bar {\bar x}}+\hat B\bar v-\dot {\bar{\bar x}}$. Then we have ${\dot {V}}(\hat e)=2\hat e^TP_0S(\varepsilon)(\hat Ae+\hat Be^v)+2\hat e^TP_0S(\varepsilon)c$. First consider $2\hat e^TP_0S(\varepsilon)(\hat Ae+\hat Be^v)$. Recall that $S(\varepsilon)=\mathrm{diag}(S_i(\varepsilon))_{i=1}^N$ and $P_0=\mathrm{diag}(P_{0i})_{i=1}^N$ are both block diagonal matrices. We then have
\begin{align}
\label{Theorem 3 step 1}
&2\hat e^TP_0S(\varepsilon)(\hat Ae+\hat Be^v)=\sum\nolimits_{i = 1}^N [2\hat e_i^T{P_{0i}}{S_i}(\varepsilon )({\hat A_{ii}} \nonumber\\
&+ \frac{1}{\varepsilon }{\hat B_i}{K_{0i}}{{\bar S}_i}(\varepsilon )){e_i}]+ \sum\nolimits_{i = 1}^N {[2\hat e_i^T{P_{0i}}{S_i}(\varepsilon )\sum\nolimits_{j \ne i} {{\hat A_{ij}}{e_j}} ]} \nonumber\\
&=\frac{2}{\varepsilon }\sum\nolimits_{i = 1}^N {[\hat e_i^T{P_{0i}}({\hat A_{ii}} + {\hat B_i}{K_{0i}}){{\hat e}_i}]}+ \sum\nolimits_{i = 1}^N [2\hat e_i^T{P_{0i}}{S_i}(\varepsilon )\nonumber\\
&\sum\nolimits_{j \ne i} {{\hat A_{ij}}{e_j}} ]=2\hat e^TP_0S(\varepsilon)\hat A_{-ii}S^{-1}(\varepsilon)\hat e-\frac{1}{\varepsilon}{\| {\hat e} \|^2}
\end{align}
where $\hat A_{-ii}=\hat A-\mathrm{diag}(\hat A_{ii})_{i=1}^N$. For $P_0S(\varepsilon )\hat A_{-ii}S^{ - 1} (\varepsilon )$, the blocks at the diagonal positions are zero matrices and the block at the $i$-th row and $j$-th column position with $i\ne j$ is $P_{0i}S_i(\varepsilon)\hat A_{ij}S_j^{-1}(\varepsilon)$. By Assumption \ref{asm: S constraint}, for each $i,j\in \mathcal{V}$, $i\ne j$, $\|P_{0i}S_i(\varepsilon)\hat A_{ij}S_j^{-1}(\varepsilon)\|_1<\sigma_i\|P_{0i}\|_1$, which implies that the maximum absolute column sum for $P_{0i}S_i(\varepsilon)\hat A_{ij}S_j^{-1}(\varepsilon)$ is upper bounded by $\sigma_i\|P_{0i}\|_1$. Each column of $P_0S(\varepsilon )\hat A_{-ii}S^{ - 1} (\varepsilon )$ is composed of columns of $N-1$ blocks $P_{0i}S_i(\varepsilon)\hat A_{ij}S_j^{-1}(\varepsilon)$ at the off-diagonal position and one column of the zero matrix at the diagonal position.
Thus, we have ${\| P_0S(\varepsilon )\hat A_{-ii}S^{ - 1} (\varepsilon ) \|_1} \le (N-1)\max_{i\in\mathcal{V}}\{\sigma_i\|P_{0i}\|_1\}$. By 5.6.P5 of \cite{Horn}, the 2-norm of a matrix with $n$ columns is upper bounded by $\sqrt{n}$ times of its 1-norm. Hence, $\| P_0S(\varepsilon )\hat A_{-ii}S^{ - 1} (\varepsilon ) \| \le \sqrt{n}{\| P_0S(\varepsilon )\hat A_{-ii}S^{ - 1} (\varepsilon ) \|_1} \le \sqrt n (N-1)\max_{i\in\mathcal{V}}\{\sigma_i\|P_{0i}\|_1\}$. By \eqref{Theorem 3 step 1}, we have $2\hat e^TP_0S(\varepsilon)(\hat Ae+\hat Be^v)\le -\frac{1}{\varepsilon}{\| {\hat e} \|^2}+2\sqrt n (N - 1)\max_{i\in\mathcal{V}}\{\sigma_i\|P_{0i}\|_1\}{\| {\hat e} \|^2}=- (\frac{1}{\varepsilon } - 2\sqrt n (N - 1)\max_{i\in\mathcal{V}}\{\sigma_i\|P_{0i}\|_1\}){\| {\hat e} \|^2}$. For any $\|\hat e\|\geq \|P_0\|\|c\|$, given $\|S(\varepsilon)\|=1$ for any $0<\varepsilon<1$, we have
\begin{align}
\label{Theorem 3 step 3}
& 2\hat e^TP_0S(\varepsilon)(\hat Ae+\hat Be^v)+2\hat e^TP_0S(\varepsilon)c\nonumber\\
&\!\leq 2\|\hat e\|\|P_0\|\|c\| - (\frac{1}{\varepsilon } - 2\sqrt n (N - 1)\max_{i\in\mathcal{V}}\{\sigma_i\|P_{0i}\|_1\}){\| {\hat e} \|^2}\nonumber\\
&\!\leq - (\frac{1}{\varepsilon } - 2\sqrt n (N - 1)\max_{i\in\mathcal{V}}\{\sigma_i\|P_{0i}\|_1\} -2){\| {\hat e} \|^2}.
\end{align}

By \eqref{Theorem 3 step 3}, by choosing $\varepsilon$ such that $\frac{1}{\varepsilon } > 2\sqrt n (N - 1)\max_{i\in\mathcal{V}}\{\sigma_i\|P_{0i}\|_1\} +2$, we have that, for any $\|\hat e\|\geq \|P_0\|\|c\|$, ${\dot {V}}(\hat e)\leq - (\frac{1}{\varepsilon } - 2\sqrt n (N - 1)\max_{i\in\mathcal{V}}\{\sigma_i\|P_{0i}\|_1\}-2){\| {\hat e} \|^2}<0$ for any $\| {\hat e} \|\neq0$. Since $P_0$ is symmetric, we have $\lambda_{\min}(P_0)\| {\hat e} \|^2\leq V(\hat e)\leq\lambda_{\max}(P_0)\| {\hat e} \|^2$, where $\lambda_{\min}(P_0)$ and $\lambda_{\max}(P_0)$ are the smallest and largest eigenvalues of $P_0$, respectively. Because $P_0$ is positive definite, $\lambda_{\min}(P_0)$ and $\lambda_{\max}(P_0)$ are both positive. View $\hat e$ as state and $c$ as input. Then, by Theorem 10.4.1 of \cite{AI:99}, the system is input-to-state stable. By Theorem 10.4.5 of \cite{AI:99}, we then have $\mathop {\lim \sup}\nolimits_{t \to \infty }\|\hat e(t)\|\leq \gamma (\mathop {\lim \sup}\nolimits_{t \to \infty }\|c(t)\|)$, where $\gamma(\cdot)$ is a class $\mathcal{K}$ function (\cite{AI:99}, Definition 10.1.1). By Theorem \ref{theorem auxiliary variable convergence}, we have that $(\bar x(t),\bar u(t))$ asymptotically converges to the VE $(\tilde x,\tilde u$) which satisfies $A\tilde x+B\tilde u=0_n$. Thus, as $t\to \infty$, we have $A\bar x(t)+B\bar u(t)\to 0_n$ and $\dot {\bar x}(t)\to 0_n$, by which we can obtain $\hat A\bar{\bar x}+\hat B\bar v\to0_n$ and $\dot{\bar{\bar x}}(t)\to0_n$. Therefore, we have $\mathop {\lim \sup}\nolimits_{t \to \infty }\|c(t)\|=0$, which implies that $\mathop {\lim \sup}\nolimits_{t \to \infty }\|\hat e(t)\|=0$. Hence, $\hat e(t)$ asymptotically converges to zero which implies that $x(t)$ keeps track of $\bar x(t)$ asymptotically. Thus, $e^v(t)$ diminishes asymptotically which implies that $v(t)$ asymptotically keeps track of $\bar v(t)$, which further implies that $u(t)$ asymptotically keeps track of $\bar u(t)$. Hence, $(x(t),u(t))$ also asymptotically converge to the VE $(\tilde x,\tilde u)$. This completes the proof of Theorem \ref{theorem: Algorithm 1}.
\end{IEEEproof}

\subsection{Proof of Lemma \ref{lemma: coordinate transformation dynamics}\label{proof of lemma}}

In this subsection, we provide the proof of Lemma \ref{lemma: coordinate transformation dynamics}.

\begin{IEEEproof}
By the construction of $z_{i1},z_{i2},z_{i3}$ in Section \ref{section nonlinear algorithm design}, the dynamics of $z_{i1}$ and $z_{i2}$ satisfy \eqref{transformed states dynamics}. Assume that the dynamics of $z_{i1},\cdots,z_{i(\ell-1)}$ satisfy \eqref{transformed states dynamics} for some $\ell\in\{3,\cdots,\bar n-1\}$. By \eqref{coordinate transformation}, we then have
\begin{align*}
&\dot z_{i\ell}=
(-1)^{\ell-1}T_{\ell-2,1}^i(z_{i3}-k_{i2}z_{i2})+(-1)^{\ell-2}T_{\ell-3,2}^i(z_{i4}-k_{i3}\nonumber\\
&z_{i3})+\cdots+(-1)^2T_{1,\ell-2}^i(z_{i\ell}-k_{i(\ell-1)}z_{i(\ell-1)})+\frac{d^{\ell-1}}{dt^{\ell-1}}\nonumber\\
&(k_1\nabla_{z_{i1}}L_i+\Gamma_{i1})+\frac{d^{\ell-2}}{dt^{\ell-2}}(x_{i2}\frac{d}{dt}\Theta_{i1}+\Theta_{i1}\Gamma_{i2})+\cdots+\nonumber\\
&\frac{d^2}{dt^2}(x_{i(\ell-2)}\frac{d}{dt}\Theta_{i}^{[\ell-3]}+\Theta_{i}^{[\ell-3]}\Gamma_{i(\ell-2)})+\frac{d}{dt}(x_{i(\ell-1)}\frac{d}{dt}\Theta_{i}^{[\ell-2]}\nonumber\\
&+\Theta_{i}^{[\ell-2]}\Gamma_{i(\ell-1)})+x_{i\ell}\frac{d}{dt}\Theta_{i}^{[\ell-1]}+\Theta_{i}^{[\ell-1]}(\Gamma_{i\ell}+\Theta_{i\ell}x_{i(\ell+1)})\nonumber\\
&=-k_{i\ell}z_{i\ell}+z_{i(\ell+1)}.
\end{align*}

Thus, the dynamics of $z_{i\ell}$ also satisfies \eqref{transformed states dynamics}. By induction, the dynamics of $z_{i\ell}$ satisfies \eqref{transformed states dynamics} for all $\ell=1,\cdots,\bar n-1$. Finally, one can expand the dynamics of $z_{i\bar n}$ as above and, by \eqref{control law for nonlinear dynamics}, obtain $\dot z_{i\bar n}=-k_{i\bar n}z_{i\bar n}$.
\end{IEEEproof}

\subsection{Proof of Theorem \ref{theorem: Algorithm 2}\label{proof of theorem 2}}

This subsection proves Theorem \ref{theorem: Algorithm 2}.

\begin{IEEEproof} First we show the convergence of the transformed system in the following theorem.

\btheorem
\label{theorem z system}
For any solution of \eqref{nonlinear state dynamics}, \eqref{coordinate transformation}, \eqref{control law for nonlinear dynamics} and \eqref{nonlinear dual update} starting from any $(x(0),u(0),\lambda(0))\in\mathbb{R}^n\times\mathbb{R}^N\times\mathbb{R}^m$, $(z,\lambda)$ asymptotically converges to a point $(\tilde z,\tilde\lambda)$ where $\tilde z^{[1]}$ is the VE of problem \eqref{OP nonlinear} and $\tilde z^{[\ell]}=0_N$ for all $\ell=2,\cdots,\bar n$.
\etheorem

\begin{IEEEproof} Let $\nabla_{z^{[1]}}L(z^{[1]},\lambda)=(\nabla_{z_{i1}}L_i(z^{[1]},\lambda))_{i=1}^N$. Under Assumptions \ref{asm: convexity}--\ref{asm: Slater condition}, by Section 12.2.3 of \cite{PE:10}, $\tilde x^{[1]}\in \mathbb{R}^N$ is the VE of the GNEP$(\mathcal{D},F)$ if and only if there exist $\tilde \lambda\in\mathbb{R}^m$ such that the KKT conditions hold: $\nabla_{z^{[1]}}L(\tilde x^{[1]},\tilde\lambda)=0_N$ and $H(\tilde x^{[1]})=0_m$.
Denote by $Q$ the set of $(\tilde x^{[1]},\tilde\lambda)$ satisfying the KKT conditions. Let $\tilde z_{i1}=\tilde x_{i1}$, $\tilde z_{i\ell}=0$ for all $\ell=2,\cdots,\bar n$, $y=[z^T,\lambda^T]^T$ and $\tilde y=[\tilde z^T,\tilde\lambda^T]^T$ with $(\tilde x^{[1]},\tilde\lambda)$ being any point of $Q$.
Choose the Lyapunov function as $V(y)=\frac{1}{2}\|y-\tilde y\|^2$. Let $\mathcal{L}V(y)$ be the Lie derivative of $V$ along \eqref{nonlinear state dynamics}, \eqref{coordinate transformation}, \eqref{control law for nonlinear dynamics} and \eqref{nonlinear dual update} at $y\in\mathbb{R}^n\times\mathbb{R}^m$. We next show that $V$ is non-increasing along any solution of \eqref{nonlinear state dynamics}, \eqref{coordinate transformation}, \eqref{control law for nonlinear dynamics} and \eqref{nonlinear dual update}. By Lemma \ref{lemma: coordinate transformation dynamics}, we have
\begin{align}
\label{dot V}
&\mathcal{L}V(y)=\sum\nolimits_{i\in\mathcal{V}}[(z_{i1}-\tilde z_{i1})\dot z_{i1}+z_{i2}\dot z_{i2}+\cdots+z_{i\bar n}\dot z_{i\bar n}]\nonumber\\
&+(\lambda-\tilde\lambda)^T\dot\lambda=\sum\nolimits_{i\in\mathcal{V}}[(z_{i1}-\tilde z_{i1})(z_{i2}-k_{1}\nabla_{z_{i1}}L_i(y))\nonumber\\
&+z_{i2}(z_{i3}-k_{i2}z_{i2})+\cdots+z_{i(\bar n-1)}(z_{i\bar n}-k_{i(\bar n-1)}z_{i(\bar n-1)})\nonumber\\
&-z_{i\bar n}k_{i\bar n}z_{i\bar n}]+(\lambda-\tilde\lambda)^Tk_1H(z^{[1]})=-\sum\nolimits_{i\in\mathcal{V}}k_1(z_{i1}-\tilde z_{i1})\nonumber\\
&\nabla_{z_{i1}}L_i(y)+k_1(\lambda-\tilde\lambda)^TH(z^{[1]})+\sum\nolimits_{i\in\mathcal{V}}[-k_{i2}z_{i2}^2-\cdots-\nonumber\\
&k_{i\bar n}z_{i\bar n}^2+(z_{i1}-\tilde z_{i1})z_{i2}+z_{i2}z_{i3}+\cdots+z_{i(\bar n-1)}z_{i\bar n}].
\end{align}
By the affinity of $H$ and the KKT condition $\nabla_{z^{[1]}}L(\tilde x^{[1]},\tilde\lambda,\tilde\mu)=0_N$, we can derive
\begin{align*}
&-\sum\nolimits_{i\in\mathcal{V}}(z_{i1}-\tilde z_{i1})\nabla_{z_{i1}}L_i(y)=-(z^{[1]}-\tilde z^{[1]})^T(\nabla F(z^{[1]})\nonumber\\
&-\nabla F(\tilde z^{[1]}))-(\lambda-\tilde\lambda)^T(H(z^{[1]})-H(\tilde z^{[1]})).
\end{align*}
Then, by strong monotonicity of $\nabla F$, we have
\begin{align}
\label{dot V z1}
&-\sum\nolimits_{i\in\mathcal{V}}k_1(z_{i1}-\tilde z_{i1})\nabla_{z_{i1}}L_i(y)+k_1(\lambda-\tilde\lambda)^TH(z^{[1]})\nonumber\\
&=-k_1(z^{[1]}-\tilde z^{[1]})^T(\nabla F(z^{[1]})-\nabla F(\tilde z^{[1]}))\nonumber\\
&\leq-k_1M\sum\nolimits_{i\in\mathcal{V}}(z_{i1}-\tilde z_{i1})^2.
\end{align}

For each $i\in\mathcal{V}$, let $\hat z_{i1}=z_{i1}-\tilde z_{i1}$. Given that $k_1,k_{i2},\cdots,k_{i\bar n}$ are chosen such that $k_1M>1,k_{i2}>1,\cdots,k_{i\bar n}>1$, plugging \eqref{dot V z1} into \eqref{dot V} yields
\begin{align}
\label{dot V overall}
&\mathcal{L}V(y)\leq-k_1M\sum\nolimits_{i\in\mathcal{V}}\hat z_{i1}^2+\sum\nolimits_{i\in\mathcal{V}}[-k_{i2}z_{i2}^2-\cdots-k_{i\bar n}z_{i\bar n}^2\nonumber\\
&+\hat z_{i1}z_{i2}+z_{i2}z_{i3}+\cdots+z_{i(\bar n-1)}z_{i\bar n}]\leq-\sum\nolimits_{i\in\mathcal{V}}[(k_1M-1)\nonumber\\
&\hat z_{i1}^2+(k_{i2}-1)z_{i2}^2+\cdots+(k_{i\bar n}-1)z_{i\bar n}^2]\leq0.
\end{align}

This proves a counterpart result of Claim \ref{middle lemma monotonicity of Lyapunov function}. Following the same arguments, one can show that counterparts of Claim \ref{middle lemma existence and uniqueness} and Claim \ref{lemma invariance omega limit set} hold for \eqref{nonlinear state dynamics}, \eqref{coordinate transformation}, \eqref{control law for nonlinear dynamics} and \eqref{nonlinear dual update}. We are ready to apply Lemma \ref{lemma invariance}. For any $\delta>0$, consider the compact set $S=V^{-1}(\leq \delta)\cap(\mathbb{R}^n\times\mathbb{R}^m)$. By Lemma \ref{lemma invariance}, starting from any point in $S$, the solution of \eqref{transformed states dynamics} and \eqref{nonlinear dual update} converges to the largest invariant set $M$ contained in ${\rm cl}(Z)$, where $Z=\{(z,\lambda)\in S:\mathcal{L}V(z,\lambda)=0\}$. By the above proof, we have $\mathcal{L}V(z,\lambda)=0$ if and only if $z=\tilde z$. Hence, $Z=\{(z,\lambda)\in S:z=\tilde z\}$. It is clear that $Z$ is closed. Let $(\tilde z,\lambda)\in M\subseteq Z$. The solution of \eqref{transformed states dynamics} and \eqref{nonlinear dual update} starting from $(\tilde z,\lambda)$ remains in $M$ only if $\nabla_{z^{[1]}}L(\tilde z^{[1]},\lambda)=0_N$. Together with $H(\tilde z^{[1]})=0_m$, we have that the pair $(\tilde z^{[1]},\lambda)$ satisfies the KKT conditions. Hence, we have $M\subseteq Q$. By the definition of omega-limit set and \eqref{dot V overall}, the omega-limit set of any solution of \eqref{transformed states dynamics} and \eqref{nonlinear dual update} is a singleton. Since $\delta$ is arbitrary, we have that, starting from any point in $\mathbb{R}^n\times\mathbb{R}^m$, the solution of \eqref{transformed states dynamics} and \eqref{nonlinear dual update} converges to a point $(\tilde z,\tilde\lambda)$ with $(\tilde z^{[1]},\tilde\lambda)$ being a point in $Q$ and $\tilde z^{[\ell]}=0_N$ for all $\ell=2,\cdots,\bar n$.
\end{IEEEproof}

We next show the convergence of the original system. Since $z_{i1}=x_{i1}$, $x^{[1]}$ asymptotically converges to the VE $\tilde x^{[1]}$. Since $\tilde z_{i2}=0$ and $\nabla_{z_{i1}}L_i(\tilde y)=0$, by the coordinate transformation of $z_{i2}$, $\Theta_{i1}(\tilde x^{[1]})\tilde x_{i2}+\Gamma_{i1}(\tilde x^{[1]})=0$. This is the steady state condition for $\dot x_{i1}$. We next show by induction that, for each $\ell=2,\cdots,\bar n$, $x_{i\ell}$ asymptotically converges to $\tilde x_{i\ell}$, the steady state of $x_{i\ell}$ given $\tilde x^{[1]}$. We have proved the case of $\ell=2$. Assume that this argument holds for all $\ell$ up some $k\in\{2,\cdots,\bar n-2\}$. For $\ell=k+1$, we check the coordinate transformation of $z_{i(\ell+1)}$ given by \eqref{coordinate transformation}. Notice that $\tilde z_{i2}=\cdots=\tilde z_{i(k+1)}=0$. Also notice that on the right-hand-side of the equation, all the terms with the derivative operation $\frac{d}{dt}$ only depend on $x^{[1\sim k]}$, which, by assumption, asymptotically converge to fixed values $\tilde x^{[1\sim k]}$. Thus, the derivatives asymptotically converge to zero. Therefore, as $t\to\infty$, $\Theta_{i}^{[\ell-1]}\Gamma_{i\ell}+\Theta_{i}^{[\ell]}x_{i(\ell+1)}\to0$, which implies $\Gamma_{i\ell}+\Theta_{i\ell}x_{i(\ell+1)}\to0$. This is the steady state condition for $\dot x_{i\ell}$. Thus, the argument holds for $\ell=k+1$. By induction, the argument holds for all $\ell=2,\cdots,\bar n-1$. Similarly, for \eqref{control law for nonlinear dynamics}, as $t\to\infty$, since $\tilde z_{i2}=\cdots=\tilde z_{i\bar n}=0$ and all the derivative terms diminish, we have $u_i\to-\frac{\Gamma_{i\bar n}}{\Theta_{i\bar n}}$, which is the steady state condition for $\dot x_{i\bar n}$. Thus, for all $i\in\mathcal{V}$, $x_{i2},\cdots,x_{i\bar n},u_i$ asymptotically converge to $\tilde x_{i2},\cdots,\tilde x_{i\bar n},\tilde u_i$ such that, together with $\tilde x^{[1]}$, $(\tilde x,\tilde u)$ is a steady state of system \eqref{nonlinear state dynamics}.
\end{IEEEproof}

\section{Case study\label{simulation section}}

In this section, we validate Theorem \ref{theorem: Algorithm 1} and Theorem \ref{theorem: Algorithm 2} by the multi-zone building temperature regulation problem and the OPF problem introduced in Section \ref{section motivation}, respectively.

\subsection{Simulation for Algorithm \ref{Algorithm 1}\label{OPF simulation}}

We simulate Algorithm \ref{Algorithm 1} by the multi-zone building temperature regulation problem introduced in Section \ref{section MZB}. Assume that the communication graph satisfies Assumption \ref{asm: general communication graph}.

\emph{System dynamics}. In \eqref{thermal dynamics}, the states of zone $i$ are $(T_{i1},T_{i2})$; the control inputs of zone $i$ are $(m_i^s,\Delta T_i^h)$; the control inputs of the AHU are $(\Delta T^c,\delta)$; the terms $(P_i^d,T^{oa})$ are disturbances which are assumed to be perfectly measured. We make the following simplifications so that the states of each zone $i$ are controlled by a single control input: 1) $\delta(t)\equiv0$, i.e., the return air is not recirculated; 2) $\Delta T^c(t)\equiv0$, i.e., the cooling power of the AHU is not used, instead, we assume that zone $i$'s control $\Delta T_i^h$ could generate both heating and cooling powers; 3) $m_i^s(t)\equiv \bar m_i^s$ where $\bar m_i^s$ is a constant known to zone $i$, i.e., the mass flow rate of the air supplied to zone $i$ is constant and known. We also assume that $P_i^d(t)\equiv\bar P_i^d$ and $T^{oa}(t)\equiv\bar T^{oa}$ where $\bar P_i^d$ and $\bar T^{oa}$ are constants known to zone $i$ and all the zones, respectively, i.e., the external load and the outside air temperature are both constant and known. For each $i\in\mathcal{V}$, let $x_{i1}=T_{i2}$, $x_{i2}=T_{i1}$, and $u_i=\bar m_i^sc_p\Delta T_i^h+\bar m_i^sc_p\bar T^{oa}+\frac{\bar T^{oa}}{R_i^{oa}}+\bar P_i^d$. Under the above simplifications, we can obtain the following dynamic model from \eqref{thermal dynamics}:
\begin{align}
\label{simplified thermal dynamics}
&\dot x_{i1}(t)=-\frac{1}{C_{i2}R_i}x_{i1}(t)+\frac{1}{C_{i2}R_i}x_{i2}(t),\nonumber\\
&\dot x_{i2}(t)=\frac{1}{C_{i1}}(\frac{x_{i1}(t)-x_{i2}(t)}{R_i}-\bar m_i^sc_px_{i2}(t)-\frac{x_{i2}(t)}{R_i^{oa}}\nonumber\\
&+\sum\nolimits_{j\in\mathcal{N}_i^{PH}}\frac{x_{j2}(t)-x_{i2}(t)}{R_{ji}})+\frac{1}{C_{i1}}u_i(t).
\end{align}

It is easy to check that \eqref{simplified thermal dynamics} satisfies Assumption \ref{asm: CCF}. We next transform \eqref{simplified thermal dynamics} into controllable canonical form. For each $i\in\mathcal{V}$, let $\mathcal{T}_i=\left[ {\begin{array}{*{20}{c}}
{C_{i1}C_{i2}R_i}&0\\
-{C_{i1}}&{C_{i1}}
\end{array}} \right]$, $\bar A_{ii}=\mathcal{T}_iA_{ii}\mathcal{T}_i^{-1}$, $\hat B_i=\mathcal{T}_iB_i$, and $\hat A_{ij}=\mathcal{T}_iA_{ij}\mathcal{T}_i^{-1}$. It can be checked that $\hat A_{ij}=\left[ {\begin{array}{*{20}{c}}
0&0\\
\frac{1}{C_{i2}R_iR_{ji}}&\frac{1}{R_{ji}}
\end{array}} \right]$ for any $j\in\mathcal{N}_i^{PH}$, and $(\bar A_{ii},\hat B_i)$ is in the controllable canonical form such that $\bar A_{ii}=[[0,1]^T,-a_i^T]^T$ with $a_i=[\frac{1}{C_{i1}C_{i2}R_i}(\bar m_i^sc_p+\frac{1}{R_i^{oa}}+\sum\nolimits_{j\in\mathcal{N}_i^{PH}}\frac{1}{R_{ji}}),\frac{1}{C_{i1}}(\frac{C_{i1}+C_{i2}}{C_{i2}R_i}+\bar m_i^sc_p+\frac{1}{R_i^{oa}}+\sum\nolimits_{j\in\mathcal{N}_i^{PH}}\frac{1}{R_{ji}})]$ and $\hat B_i=[0,1]^T$. For any $j\in\mathcal{N}_i^{PH}$, given any $0<\varepsilon<1$, it can be checked that $S_i(\varepsilon)\hat A_{ij}S_j^{-1}(\varepsilon)=\left[ {\begin{array}{*{20}{c}}
0&0\\
\frac{\varepsilon^2}{C_{i2}R_iR_{ji}}&\frac{\varepsilon}{R_{ji}}
\end{array}} \right]$. Hence, for any $0<\varepsilon<1$, we have $\|S_i(\varepsilon)\hat A_{ij}S_j^{-1}(\varepsilon)\|_1=\max\{\frac{\varepsilon^2}{C_{i2}R_iR_{ji}},\frac{\varepsilon}{R_{ji}}\}<\max\{\frac{1}{C_{i2}R_iR_{ji}},\frac{1}{R_{ji}}\}<\frac{1}{C_{i2}R_iR_{ji}}+\frac{1}{R_{ji}}=\frac{1+C_{i2}R_i}{C_{i2}R_iR_{ji}}\leq\frac{1+C_{i2}R_i}{C_{i2}R_i\min_{\ell\in\mathcal{N}_i^{PH}}\{R_{\ell i}\}}$, which is a constant independent of $j$. Hence, Assumption \ref{asm: S constraint} is satisfied.


\emph{GNEP}. We formulate \eqref{thermal OP} into a GNEP. Let $\bar c_{i1}^x=c_{i2}^x$, $\bar c_{i2}^x=c_{i1}^x$, $\bar c_i^u=\frac{c_i^u}{(\bar m_i^sc_p)^2}$, $u_i^r=\bar m_i^sc_p\bar T^{oa}+\frac{\bar T^{oa}}{R_i^{oa}}+\bar P_i^d$, $\underline u_i=\bar m_i^sc_p\underline{\Delta T}_i^h+\bar m_i^sc_p\bar T^{oa}+\frac{\bar T^{oa}}{R_i^{oa}}+\bar P_i^d$, and $\overline u_i=\bar m_i^sc_p\overline{\Delta T}_i^h+\bar m_i^sc_p\bar T^{oa}+\frac{\bar T^{oa}}{R_i^{oa}}+\bar P_i^d$. Given $(x_{-i},u_{-i})$, agent $i$ solves the following optimization problem
\begin{align}
\label{thermal OP simulation with steady state conditions}
&\!\!\!\!\!\min\limits_{x_{i1},x_{i2}\in[\underline T_i,\overline T_i],u_i\in[\underline u_i,\overline u_i]}\bar c_{i1}^x(x_{i1}-T_i^{r})^2+\bar c_{i2}^x(x_{i2}-T_i^{r})^2\nonumber\\
&\!\!\!\!\!\quad\quad\quad\quad\quad\quad\quad\quad\quad\;\,+\bar c_{i}^u(u_i-u_i^r)^2\nonumber\\
&\!\!\!\!\!{\rm s.t.}\;\; u_i-\bar m_i^sc_px_{i2}-\frac{x_{i2}}{R_i^{oa}}+\sum\nolimits_{j\in\mathcal{N}_i^{PH}}\frac{x_{j2}-x_{i2}}{R_{ji}}=0,\nonumber\\
&\!\!\!\!\!\quad\;\;\; x_{i1}-x_{i2}=0.
\end{align}
It is easy to check that problem \eqref{thermal OP simulation with steady state conditions} satisfies Assumptions \ref{asm: convexity}--\ref{asm: strict monotone} and Assumptions \ref{asm: W bounded}--\ref{asm: physical feasibility}.
In the above, we have verified that all the assumptions needed by Theorem \ref{theorem: Algorithm 1} are satisfied.

\begin{figure}
\begin{center}
\includegraphics[width=1\linewidth]{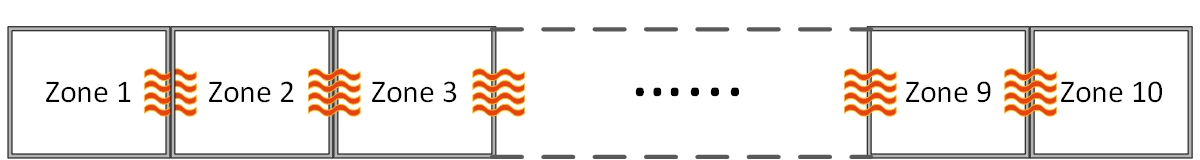}
\caption{Floor plan for the temperature regulation simulation}
\label{floor_plan}
\end{center}
\end{figure}

\begin{figure}
\begin{center}
\includegraphics[width=1\linewidth]{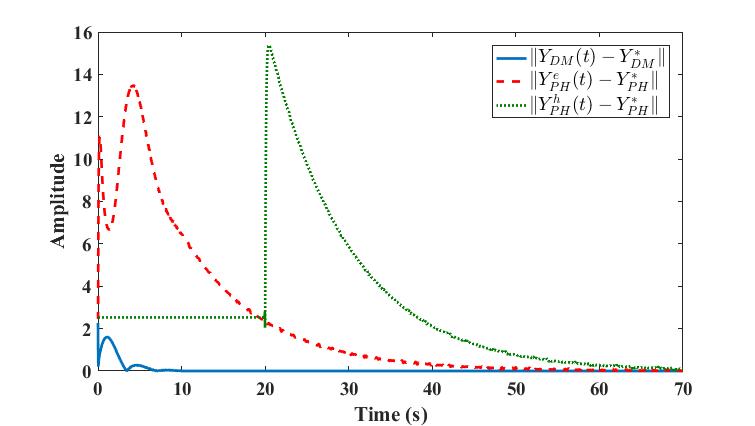}
\caption{2-norm error of the auxiliary variables and the true variables}
\label{DM_PH_err}
\end{center}
\end{figure}

\emph{Simulation results}. We consider the case where $N=10$. The undirected graph describing the physical topology of the zone network is denoted by $\mathcal{G}^{PH}=(\mathcal{V},\mathcal{E}^{PH})$, where $\mathcal{V}=\{1,2,\cdots,10\}$ and $\mathcal{E}^{PH}=\{(1,2),(2,3),\cdots,(8,9),(9,10)\}$. The floor plan is depicted by Fig. \ref{floor_plan}. This adjacency topology is widely used for case studies in the literature, e.g., \cite{YM-GA-FB:2011}. The values of the parameters are adopted from \cite{YM-GA-FB:2011}: for all $i\in\mathcal{V}$, $C_{i1}=9163kJ/K$, $C_{i2}=169400kJ/K$, $m_i^s=0.01kg/s$, $T^{oa}=25\degree C$, $R_i=1.7K/kW$, $R_{ji}=2K/kW$, $R_i^{oa}=57K/kW$, $c_p=1012J/kg\cdot K$, $T_i^{ref}=21.6\degree C$, $\underline T_i=20.6\degree C$, $\overline T_i=21.7\degree C$, $\underline{\Delta T}_i^h=-30\degree C$, $\overline{\Delta T}_i^h=8\degree C$ and $P_i^d=0.1kW$. We choose $c_{i1}^x=c_{i2}^x=10$, and $c_{i}^u=0.1$. Let $(\bar x,\bar u)$ be the auxiliary variables, $Y_{DM}=[\bar x_{11},\bar x_{12},\bar u_1,\cdots,\bar x_{N1},\bar x_{N2},\bar u_N]^T$ and $Y_{DM}^*$ be the VE of problem \eqref{thermal OP simulation with steady state conditions}. For the true variables $(x,u)$, let $Y_{PH}=[x_{11},x_{12},u_1,\cdots,x_{N1},x_{N2},u_N]^T$ and $Y_{PH}^*$ be the steady state of $Y_{PH}$. We use $Y_{PH}^e(t)$ and $Y_{PH}^h(t)$ to denote the trajectories of $Y_{PH}$ obtained by our economic control scheme and the hierarchical control scheme, respectively. In Fig. \ref{DM_PH_err}, the 2-norm errors $\|Y_{DM}(t)-Y_{DM}^*\|$, $\|Y_{PH}^e(t)-Y_{PH}^*\|$ and $\|Y_{PH}^h(t)-Y_{PH}^*\|$ are shown by the blue solid line, the red dashed line and the green dot line, respectively. This verifies the convergence of the algorithm. The areas under the trajectories of $\|Y_{PH}^e(t)-Y_{PH}^*\|$ and $\|Y_{PH}^h(t)-Y_{PH}^*\|$, denoted by $S_e$ and $S_h$, represent efficiency loss of economic control and hierarchical control in reaching convergence, respectively. In the simulation, we have $S_e=158.5730$ and $S_h=207.6105$. This verifies that the proposed economic control scheme can significantly reduce efficiency loss.

\subsection{Simulation for Algorithm \ref{Algorithm 2}\label{section simulation problem OPF}}

\begin{figure}
\begin{center}
\includegraphics[width=1.8in]{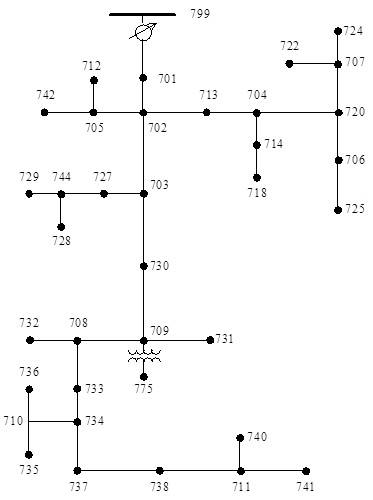}
\caption{IEEE 37-bus Test System}
\label{IEEE37}
\end{center}
\end{figure}

We simulate Algorithm \ref{Algorithm 2} by the OPF problem in Section \ref{section OPF} and use the IEEE 37-bus Test System \cite{IEEE-1991} for the agents' physical couplings; see Fig. \ref{IEEE37}. Assume that the communication graph satisfies Assumption \ref{asm: nonlinear communication graph}. The parameter values are adopted from \cite{YW-DJH-GG:1998}: $D_i=5MW/Hz$, $T_{H_i}=4s$, $T_{M_i}=0.35$, $K_{M_i}=1$, $K_{E_i}=1$, $T_{E_i}=0.1s$, $R_i=0.05Hz/MW$, $t_{ij}=1.5MW/rad$,
and $\tilde\omega_i=50deg/s$. We choose $a_i=a=1$, $b_i=b=10$ and $c_i=c=50$ for all $i\in\mathcal{V}$.


\emph{System dynamics}. In \eqref{motivating example dynamics}, agent $i$'s state is $x_i=[\theta_i,\omega_i,P_{m_i},P_{E_i}]^T$ and control input is $u_i$. It is clear that \eqref{motivating example dynamics} is in the form of \eqref{nonlinear state dynamics} and satisfies Assumptions \ref{asm: smoothness} and \ref{asm: convex feasible set}.

\emph{GNEP.} We formulate \eqref{OPF formulation} into a GNEP. Given $(P_{m_{-i}},\theta_{-i})$, agent $i$ solves the following optimization problem
\begin{align}
\label{OPF formulation GNEP}
&\mathop {\min }\nolimits_{(P_{m_i},\theta_i) \in {\mathbb{R}^{2}}} {aP_{m_i}^2+bP_{m_i}+c\theta_i^2}\nonumber\\
&{\rm{s.t.}}\sum\nolimits_{j\in\mathcal{N}_i^{PH}}{t_{ij}}(\theta_{ij}-\sin\theta_{ij}(0))=P_{m_i}-P_{m_i0}-\frac{D_i\tilde\omega_i}{\omega_0}.
\end{align}

We next transform the GNEP \eqref{OPF formulation GNEP} into the one that only depends on $x^{[1]}=(\theta_{j1})_{j=1}^N$. We use the steady state condition of $\omega_i(t)$ in \eqref{motivating example dynamics} to obtain the manifold of $P_{m_i}$ as a function of $(\theta_{j1})_{j=1}^N$ and substitute the manifold into \eqref{OPF formulation GNEP}. To maintain convexity in the game formulation, we use the linearization $\sin\theta_{ij}\approx \theta_{ij}$. Let $o_i=D_i\frac{\tilde\omega_i}{\omega_0}+P_{m_i0}-\sum\nolimits_{j\in\mathcal{N}_i^{PH}}{t_{ij}}\sin\theta_{ij}(0)$. By the dynamics of $\omega_i(t)$ in \eqref{motivating example dynamics} with linearization, the manifold of $P_{m_i}$ is $P_{m_i}=o_i+\sum\nolimits_{j\in\mathcal{N}_i^{PH}}{t_{ij}}\theta_{ij}$. Notice that the first constraint of \eqref{OPF formulation GNEP} coincides with this manifold and thus can be removed. Replacing $P_{m_i}$ in \eqref{OPF formulation GNEP} by its manifold yields
\begin{align}
\label{OPF GNEP simulation}
\mathop {\min }\limits_{\theta_i\in\mathbb{R}} {a(o_i+\sum\limits_{j\in\mathcal{N}_i^{PH}}{t_{ij}}\theta_{ij})^2+b(o_i+\sum\limits_{j\in\mathcal{N}_i^{PH}}{t_{ij}}\theta_{ij})}+c\theta_i^2.
\end{align}

It is clear that the GNEP \eqref{OPF GNEP simulation} satisfies Assumptions \ref{asm: convexity}, \ref{asm: Slater condition}, \ref{asm: Lipschitz continuity} and \ref{asm: smoothness}.
We next verify Assumption \ref{asm: strong monotone}. Denote by $f_i$ the objective function of \eqref{OPF GNEP simulation} and $\nabla F=(\nabla_{\theta_i}f_i)_{i=1}^N$. Let $t_i=\sum_{j\in\mathcal{N}_i^{PH}}t_{ij}$.
For any $\theta,\theta'\in\mathbb{R}^N$, we then have
\begin{align*}
&(\theta-\theta')^T(\nabla F(\theta)-\nabla F(\theta'))=2a\sum\nolimits_{i\in\mathcal{V}}[t_i^2(\theta_i-\theta_i')^2\nonumber\\
&-t_i\sum\nolimits_{j\in\mathcal{N}_i^{PH}}t_{ij}(\theta_i-\theta_i')(\theta_j-\theta_j')]+2c\sum\nolimits_{i\in\mathcal{V}}(\theta_i-\theta_i')^2\nonumber\\
&=a\sum\nolimits_{i\in\mathcal{V}}\sum\nolimits_{j\in\mathcal{N}_i^{PH}}t_it_{ij}[(\theta_i-\theta_i')-(\theta_j-\theta_j')]^2\nonumber\\
&+\sum\nolimits_{i\in\mathcal{V}}(2c+a(t_i^2-\sum\nolimits_{j\in\mathcal{N}_i^{PH}}t_jt_{ji}))(\theta_i-\theta_i')^2\nonumber\\
&\geq\min\nolimits_{i\in\mathcal{V}}(2c+a(t_i^2-\sum\nolimits_{j\in\mathcal{N}_i^{PH}}t_jt_{ji}))\|\theta-\theta'\|^2.
\end{align*}
By the choice of $a$, $b$, $c$ and $t_{ij}$, it can be easily checked that $2c+a(t_i^2-\sum_{j\in\mathcal{N}_i^{PH}}t_jt_{ji})>0$ for all $i\in\mathcal{V}$. By Definition \ref{def: strict monotone}, $\nabla F$ is strongly monotone on $\mathbb{R}^N$. In the above, we have verified that all the assumptions of Theorem \ref{theorem: Algorithm 2} are satisfied.

\begin{figure}[ht]
\begin{center}
\includegraphics[width=1\linewidth]{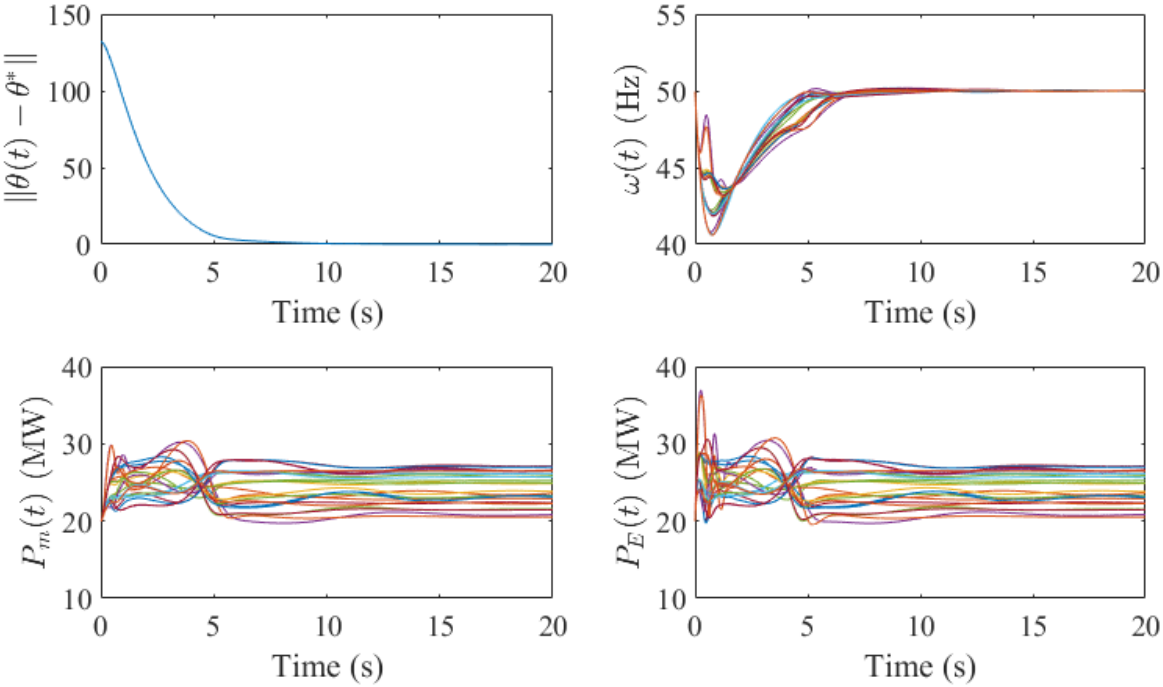}
\caption{Trajectories of the state variables}
\label{temp_regulation_err}
\end{center}
\end{figure}

\emph{Simulation results.} We choose $P_{m_i0}=20MW$, $P_{E_i0}=20MW$, $\theta_{i0}=20rad$, $\omega_{i0}=49.95 Hz$, $k_1=0.1$, and $k_{i2}=k_{i3}=k_{i4}=1.2$ for all $i\in\mathcal{V}$. Let $\theta^*$ be the VE of problem \eqref{OPF GNEP simulation}. The trajectories of $\|\theta(t)-\theta^*\|$, $\omega(t)$, $P_{m}(t)$ and $P_{E}(t)$ are given in Fig. \ref{temp_regulation_err}. This verifies that the physical variables converge to a set of steady states of system \eqref{motivating example dynamics} and the converging steady state of $\theta$ is the VE of problem \eqref{OPF GNEP simulation}.

\section{Conclusion}

We study two complementary cases of distributed economic control problems: linear systems with controllable subsystems and nonlinear systems in the strict-feedback form. The proposed algorithms are verified by two case studies on a multi-zone building temperature regulation problem and an optimal power flow problem, respectively.

\bibliographystyle{plain}
\bibliography{MZ,alias,efmain,frazzoli,YL}
\end{document}